\documentclass{amsart}
\usepackage{amsaddr}
\usepackage{orcidlink}
\usepackage{graphicx} 
\usepackage{amsfonts,amsmath,amsthm,amssymb}
\usepackage{enumerate}
\usepackage{hyperref}
\usepackage{nicefrac}
\usepackage{tikz}
\usetikzlibrary{cd}
\usepackage{caption}
\usepackage{subcaption}
\usepackage{wrapfig}

\theoremstyle{definition}
\newtheorem{theorem}{Theorem}[section]
\newtheorem{definition}[theorem]{Definition}
\newtheorem{proposition}[theorem]{Proposition}

\newtheorem{lemma}[theorem]{Lemma}
\theoremstyle{remark}
\newtheorem{remark}[theorem]{Remark}

\newcommand\sign{\mathrm{sign}}
\newcommand\conv{\mathrm{conv}}
\newcommand{\rotvert}{\rotatebox[origin=c]{90}{$\vert$}}

\hypersetup{
    colorlinks=true,
    urlcolor=blue,
    linkcolor=black
}

\usepackage{myDiagrams}
\tikzstyle{vertex}=[fill, circle,inner sep=0pt,minimum size=5pt]
\tikzstyle{vertexSmall}=[fill, circle,inner sep=0pt,minimum size=3pt]
\tikzstyle{vertexLabeled}=[draw, thick, circle,minimum size=5pt]

\begin{document}

\title{Star Triangulations}

\author{Dimitris Bogiokas\orcidlink{0000-0001-8920-5544}\and Franziska Donhauser}
\address{Department of Informatics\\TU Berlin, Germany}
\email{d.bogiokas@tu-berlin.de}
\email{franzidonhauser@googlemail.com}
\thanks{Funded, in part, by the Deutsche Forschungsgemeinschaft (DFG, German Research Foundation) under Germany´s Excellence Strategy – The Berlin Mathematics Research Center MATH+ (EXC-2046/1, project ID: 390685689); and the European Research Council (ERC) under the European Union’s Horizon 2020 research and innovation programme (Grant agreement No. 101055448, ERC Advanced Grand EMERGE).}

\keywords{Triangulations, Enumerative Combinatorics, Point Configurations, Oriented Matroids, Catalan Numbers}

\begin{abstract}
We study full triangulations of planar point configurations consisting of $n$ points in convex position and one point in the interior, assuming that no three points are collinear (star configurations). Our main contribution is a classification of the star configurations up to combinatorial equivalence. We show that there are $2^{n-1}-n$ equivalence classes, by establishing a bijection between a star configuration and a particular Young diagram. By transforming each such Young diagram to a binary sequence, we manage to also give a combinatorial description of the orbits of these star configurations under the dihedral group action, encoding each such orbit through a black-and-white bracelet.

For each star configuration, we derive an explicit formula for the number of its full triangulations, expressed in terms of Catalan numbers and the parts of the associated Young diagram. Additionally we derive a formula of these full triangulations expressed only in terms of Catalan numbers and the bead structure of the associated bracelet.

Finally, we let the points move freely under the restriction that they stay in a star configuration and show that there exist $\binom{2n-3}{n-3}$ topologically distinct full triangulations.
\end{abstract}
\maketitle
\section{Introduction}
This paper concerns the enumeration of all full triangulations of a planar point configuration $A$ consisting of $n$ boundary points and one interior point. We assume that these $n+1$ points are in general position, in the sense that no three points in $A$ are collinear.

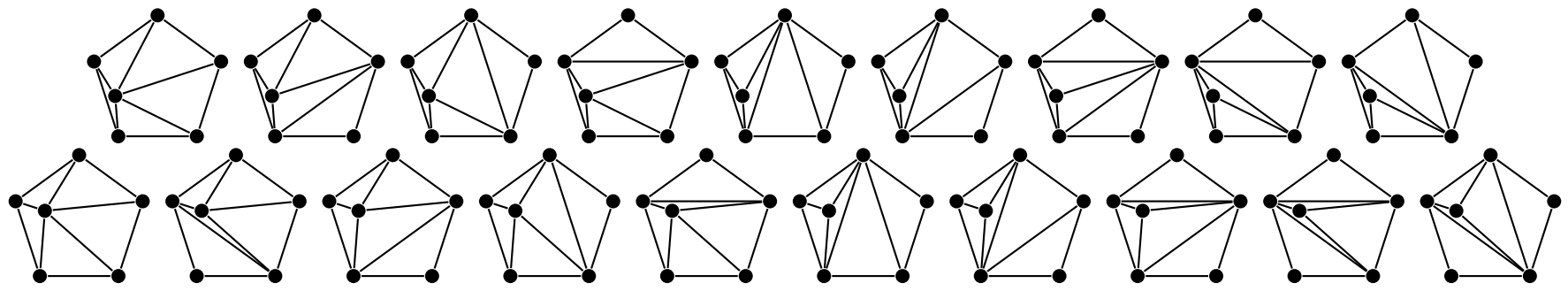
\begin{figure}[ht]
\centering
    \begin{tikzpicture}
    \def\r{0.5}
    \foreach\i in {1,2,3,4,5}
        \node[vertexSmall] (v\i) at (3*72-\i*72+18:\r) {};
    \draw (v5)--(v1)--(v2)--(v3)--(v4)--(v5);
    \def\ro{0.66666666666} 
    \node[vertexSmall] (v) at (18+2*72+36:\ro*\r) {};
    \draw (v)--(v1);
    \draw (v)--(v2);
    \draw (v)--(v3);
    \draw (v)--(v4);
    \draw (v)--(v5);
    \end{tikzpicture}
    \begin{tikzpicture}
    \def\r{0.5}
    \foreach\i in {1,2,3,4,5}
        \node[vertexSmall] (v\i) at (3*72-\i*72+18:\r) {};
    \draw (v5)--(v1)--(v2)--(v3)--(v4)--(v5);
    \def\ro{0.66666666666} 
    \node[vertexSmall] (v) at (18+2*72+36:\ro*\r) {};
    \draw (v)--(v1);
    \draw (v)--(v2);
    \draw (v)--(v3);
    \draw (v)--(v5);
    \draw (v3)--(v5);
    \end{tikzpicture}
    \begin{tikzpicture}
    \def\r{0.5}
    \foreach\i in {1,2,3,4,5}
        \node[vertexSmall] (v\i) at (3*72-\i*72+18:\r) {};
    \draw (v5)--(v1)--(v2)--(v3)--(v4)--(v5);
    \def\ro{0.66666666666} 
    \node[vertexSmall] (v) at (18+2*72+36:\ro*\r) {};
    \draw (v)--(v1);
    \draw (v)--(v2);
    \draw (v)--(v4);
    \draw (v)--(v5);
    \draw (v2)--(v4);
    \end{tikzpicture}
    \begin{tikzpicture}
    \def\r{0.5}
    \foreach\i in {1,2,3,4,5}
        \node[vertexSmall] (v\i) at (3*72-\i*72+18:\r) {};
    \draw (v5)--(v1)--(v2)--(v3)--(v4)--(v5);
    \def\ro{0.66666666666} 
    \node[vertexSmall] (v) at (18+2*72+36:\ro*\r) {};
    \draw (v)--(v1);
    \draw (v)--(v3);
    \draw (v)--(v4);
    \draw (v)--(v5);
    \draw (v1)--(v3);
    \end{tikzpicture}
    \begin{tikzpicture}
    \def\r{0.5}
    \foreach\i in {1,2,3,4,5}
        \node[vertexSmall] (v\i) at (3*72-\i*72+18:\r) {};
    \draw (v5)--(v1)--(v2)--(v3)--(v4)--(v5);
    \def\ro{0.66666666666} 
    \node[vertexSmall] (v) at (18+2*72+36:\ro*\r) {};
    \draw (v)--(v1);
    \draw (v)--(v2);
    \draw (v)--(v5);
    \draw (v2)--(v5);
    \draw (v2)--(v4);
    \end{tikzpicture}
    \begin{tikzpicture}
    \def\r{0.5}
    \foreach\i in {1,2,3,4,5}
        \node[vertexSmall] (v\i) at (3*72-\i*72+18:\r) {};
    \draw (v5)--(v1)--(v2)--(v3)--(v4)--(v5);
    \def\ro{0.66666666666} 
    \node[vertexSmall] (v) at (18+2*72+36:\ro*\r) {};
    \draw (v)--(v1);
    \draw (v)--(v2);
    \draw (v)--(v5);
    \draw (v2)--(v5);
    \draw (v3)--(v5);
    \end{tikzpicture}
    \begin{tikzpicture}
    \def\r{0.5}
    \foreach\i in {1,2,3,4,5}
        \node[vertexSmall] (v\i) at (3*72-\i*72+18:\r) {};
    \draw (v5)--(v1)--(v2)--(v3)--(v4)--(v5);
    \def\ro{0.66666666666} 
    \node[vertexSmall] (v) at (18+2*72+36:\ro*\r) {};
    \draw (v)--(v1);
    \draw (v)--(v3);
    \draw (v)--(v5);
    \draw (v1)--(v3);
    \draw (v3)--(v5);
    \end{tikzpicture}
    \begin{tikzpicture}
    \def\r{0.5}
    \foreach\i in {1,2,3,4,5}
        \node[vertexSmall] (v\i) at (3*72-\i*72+18:\r) {};
    \draw (v5)--(v1)--(v2)--(v3)--(v4)--(v5);
    \def\ro{0.66666666666} 
    \node[vertexSmall] (v) at (18+2*72+36:\ro*\r) {};
    \draw (v)--(v1);
    \draw (v)--(v4);
    \draw (v)--(v5);
    \draw (v1)--(v3);
    \draw (v1)--(v4);
    \end{tikzpicture}
    \begin{tikzpicture}
    \def\r{0.5}
    \foreach\i in {1,2,3,4,5}
        \node[vertexSmall] (v\i) at (3*72-\i*72+18:\r) {};
    \draw (v5)--(v1)--(v2)--(v3)--(v4)--(v5);
    \def\ro{0.66666666666} 
    \node[vertexSmall] (v) at (18+2*72+36:\ro*\r) {};
    \draw (v)--(v1);
    \draw (v)--(v4);
    \draw (v)--(v5);
    \draw (v2)--(v4);
    \draw (v1)--(v4);
    \end{tikzpicture}
\\
\begin{tikzpicture}
    \def\r{0.5}
    \foreach\i in {1,2,3,4,5}
        \node[vertexSmall] (v\i) at (3*72-\i*72+18:\r) {};
    \draw (v5)--(v1)--(v2)--(v3)--(v4)--(v5);
    \def\ri{0.53934466291} 
    \node[vertexSmall] (v) at (18+2*72:\ri*\r) {};
    \draw (v)--(v1);
    \draw (v)--(v2);
    \draw (v)--(v3);
    \draw (v)--(v4);
    \draw (v)--(v5);
\end{tikzpicture}
\begin{tikzpicture}
    \def\r{0.5}
    \foreach\i in {1,2,3,4,5}
        \node[vertexSmall] (v\i) at (3*72-\i*72+18:\r) {};
    \draw (v5)--(v1)--(v2)--(v3)--(v4)--(v5);
    \def\ri{0.53934466291} 
    \node[vertexSmall] (v) at (18+2*72:\ri*\r) {};
    \draw (v)--(v1);
    \draw (v)--(v2);
    \draw (v)--(v3);
    \draw (v)--(v4);
    \draw (v1)--(v4);
\end{tikzpicture}
\begin{tikzpicture}
    \def\r{0.5}
    \foreach\i in {1,2,3,4,5}
        \node[vertexSmall] (v\i) at (3*72-\i*72+18:\r) {};
    \draw (v5)--(v1)--(v2)--(v3)--(v4)--(v5);
    \def\ri{0.53934466291} 
    \node[vertexSmall] (v) at (18+2*72:\ri*\r) {};
    \draw (v)--(v1);
    \draw (v)--(v2);
    \draw (v)--(v3);
    \draw (v)--(v5);
    \draw (v3)--(v5);
\end{tikzpicture}
\begin{tikzpicture}
    \def\r{0.5}
    \foreach\i in {1,2,3,4,5}
        \node[vertexSmall] (v\i) at (3*72-\i*72+18:\r) {};
    \draw (v5)--(v1)--(v2)--(v3)--(v4)--(v5);
    \def\ri{0.53934466291} 
    \node[vertexSmall] (v) at (18+2*72:\ri*\r) {};
    \draw (v)--(v1);
    \draw (v)--(v2);
    \draw (v)--(v4);
    \draw (v)--(v5);
    \draw (v2)--(v4);
\end{tikzpicture}
\begin{tikzpicture}
    \def\r{0.5}
    \foreach\i in {1,2,3,4,5}
        \node[vertexSmall] (v\i) at (3*72-\i*72+18:\r) {};
    \draw (v5)--(v1)--(v2)--(v3)--(v4)--(v5);
    \def\ri{0.53934466291} 
    \node[vertexSmall] (v) at (18+2*72:\ri*\r) {};
    \draw (v)--(v1);
    \draw (v)--(v3);
    \draw (v)--(v4);
    \draw (v)--(v5);
    \draw (v1)--(v3);
\end{tikzpicture}
\begin{tikzpicture}
    \def\r{0.5}
    \foreach\i in {1,2,3,4,5}
        \node[vertexSmall] (v\i) at (3*72-\i*72+18:\r) {};
    \draw (v5)--(v1)--(v2)--(v3)--(v4)--(v5);
    \def\ri{0.53934466291} 
    \node[vertexSmall] (v) at (18+2*72:\ri*\r) {};
    \draw (v)--(v1);
    \draw (v)--(v2);
    \draw (v)--(v5);
    \draw (v2)--(v5);
    \draw (v2)--(v4);
\end{tikzpicture}
\begin{tikzpicture}
    \def\r{0.5}
    \foreach\i in {1,2,3,4,5}
        \node[vertexSmall] (v\i) at (3*72-\i*72+18:\r) {};
    \draw (v5)--(v1)--(v2)--(v3)--(v4)--(v5);
    \def\ri{0.53934466291} 
    \node[vertexSmall] (v) at (18+2*72:\ri*\r) {};
    \draw (v)--(v1);
    \draw (v)--(v2);
    \draw (v)--(v5);
    \draw (v2)--(v5);
    \draw (v3)--(v5);
\end{tikzpicture}
\begin{tikzpicture}
    \def\r{0.5}
    \foreach\i in {1,2,3,4,5}
        \node[vertexSmall] (v\i) at (3*72-\i*72+18:\r) {};
    \draw (v5)--(v1)--(v2)--(v3)--(v4)--(v5);
    \def\ri{0.53934466291} 
    \node[vertexSmall] (v) at (18+2*72:\ri*\r) {};
    \draw (v)--(v1);
    \draw (v)--(v3);
    \draw (v)--(v5);
    \draw (v1)--(v3);
    \draw (v3)--(v5);
\end{tikzpicture}
\begin{tikzpicture}
    \def\r{0.5}
    \foreach\i in {1,2,3,4,5}
        \node[vertexSmall] (v\i) at (3*72-\i*72+18:\r) {};
    \draw (v5)--(v1)--(v2)--(v3)--(v4)--(v5);
    \def\ri{0.53934466291} 
    \node[vertexSmall] (v) at (18+2*72:\ri*\r) {};
    \draw (v)--(v1);
    \draw (v)--(v3);
    \draw (v)--(v4);
    \draw (v1)--(v4);
    \draw (v1)--(v3);
\end{tikzpicture}
\begin{tikzpicture}
    \def\r{0.5}
    \foreach\i in {1,2,3,4,5}
        \node[vertexSmall] (v\i) at (3*72-\i*72+18:\r) {};
    \draw (v5)--(v1)--(v2)--(v3)--(v4)--(v5);
    \def\ri{0.53934466291} 
    \node[vertexSmall] (v) at (18+2*72:\ri*\r) {};
    \draw (v)--(v1);
    \draw (v)--(v2);
    \draw (v)--(v4);
    \draw (v1)--(v4);
    \draw (v2)--(v4);
\end{tikzpicture}
\caption{The two rows illustrate two different possible placements of the interior vertex in a pentagon. The corresponding configurations have numbers of full triangulations that differ by one.}
\label{fig:different_triangulations}
\end{figure}

In contrast to the case in which all points are in convex position, the set of triangulations in this setting depends on the actual embedding of the points; see, for example, Figure~\ref{fig:different_triangulations}. Thus, in order to enumerate the triangulations, we first give a combinatorial description of all possible configurations of this type for a given number of points.

\begin{definition}
A $d$-dimensional \emph{point configuration} $A$ is a finite set of labeled points $\{a_1,\ldots,a_n\}\subseteq\mathbb{R}^d$, i.e. an injective mapping $[n]\to\mathbb{R}^d$ ($i\mapsto a_i$). We write $|A|=n$.
\end{definition}
Throughout, we follow the notation of \cite{triangulations}. In the more general setting considered in that book, it is important to regard point configurations in $\mathbb{R}^d$ as collections (multisets) of points. In our setting, however, it suffices to consider only distinct points.
\begin{definition}
    Let $A\subseteq\mathbb{R}^d$ be a $d$-dimensional point configuration, with $|A|=n$. A subset $T\subseteq\binom{[n]}{d+1}$ is called a \emph{triangulation} if it satisfies the following union and intersection properties:
    \begin{itemize}
        \item\hfill
            \(\displaystyle\bigcup_{\sigma\in T}\conv(A|_{\sigma}) = \conv(A)\)\hfill(UP)
        \item\hfill
            \(\displaystyle\forall\sigma_1,\sigma_2\in T\quad\conv(A|_{\sigma_1})\cap\conv(A|_{\sigma_2}) = \conv(A|_{\sigma_1\cap\sigma_2})\)\hfill(IP)
    \end{itemize}
    where $A|_I=\{a_i:i\in I\}$.
    The triangulation is called \emph{full} if, in addition, it satisfies:
    \begin{itemize}
        \item\hfill\(\displaystyle\bigcup_{\sigma\in T}\sigma = [n]\)\hfill(FP)
    \end{itemize}
    Let $\mathcal{T}^{\text{full}}_A$ be the set of all full triangulations of $A$.
\end{definition}

By Radon's theorem, triangulations of $A$ can be characterized combinatorially in terms of its \emph{oriented circuits} $(Z^+,Z^-)$ of $A$. In our setting  a circuit consists of $Z^+,Z^-$, two disjoint subsets of $A$ with $|Z^+|+|Z^-|=d+2$ and $\conv(Z^+)\cap\conv(Z^-)\neq\emptyset$. For a detailed definition of oriented circuits in the general case we refer the reader to \cite{triangulations}.

The set of oriented circuits defines a so-called \emph{oriented matroid}, a construction more general than the point configurations. For an introduction to oriented matroids, see \cite{orientedMatroids}. For the combinatorial characterization of triangulations using the oriented circuits of $A$ see \cite{cyclicAndHigherBruhat}.
In this paper, we use the fact that the oriented matroid associated with a point configuration $A$ is fully determined by its chirotope (see \cite{computingTrUsingOrM} and \cite{polytopeOfAllTr}).

\begin{definition}
    Let $A=\{a_1,\ldots,a_n\}\subseteq\mathbb{R}^d$ be a point configuration. The \emph{chirotope} of $A$ is the function:
    \begin{center}
    \begin{tikzcd}
        
        \chi_{A}:&[-3em]\binom{[n]}{d+1}\ar[r]&\{-1,0,+1\}\\[-2em]
        &\{i_1<\cdots<i_{d+1}\}\ar[r,mapsto]&
        \sign\det\begingroup\setlength\arraycolsep{1pt}\begin{pmatrix}
            1&\rotvert a_{i_1}\rotvert\\
            \vdots&\vdots\\
            1&\rotvert a_{i_{d+1}}\rotvert
        \end{pmatrix}\endgroup
    \end{tikzcd}
    \end{center}
\end{definition}
We usually extend $\chi_A$ to an antisymmetric function on $[n]^{d+1}$. The extended function vanishes whenever an index is repeated and otherwise changes sign according to the permutation of its arguments $\chi_{A}(i_{\pi(1)},\ldots,i_{\pi(d+1)})=\sign(\pi)\cdot\chi(i_1<\cdots<i_{d+1})$. This is equivalent to using the matrix definition on any input $(d+1)$-tuple.

Our contributions are the following:
\begin{itemize}
    \item In Section~\ref{ssec:bijection} we establish a bijection between all combinatorially different point configurations and a particular set of Young diagrams.
    \item In Theorem~\ref{thm:regions} we count these different configurations.
    \item In Theorem~\ref{thm:enumTriangFromYoung} we give an explicit formula for the number of full triangulations for every configuration, with respect to its associated Young diagram.
    \item In Theorem~\ref{thm:enumTriangFromBra} we give an explicit formula for the number of full triangulations for every class of configurations, under the action of the dihedral group, with respect to its associated bracelet, the combinatorial object that naturally describes these orbits.
    \item In Theorem~\ref{thm:enumTriangTop} we count the number of every possible full topological triangulation, under the restriction that $n$ points are on the boundary and one in the interior.
\end{itemize}
\section{Enumerating the Configurations}
We first restrict the possible labelings of the configurations. Indeed, any unlabeled configuration can be labeled so that its boundary points are numbered counterclockwise from $1$ to $n$, while its unique interior point receives the label $n+1$. Consequently, there is no need to consider all possible permutations of the boundary labels.
\begin{definition}
    An $n$-\emph{star} configuration $A=\{a_1,\ldots,a_n,a_{n+1}\}$ consists of points $a_1,\ldots,a_n$ in convex position given in counterclockwise order and an inner point $a_{n+1}$, such that no three points are collinear.
\end{definition}
\subsection{\texorpdfstring{$n$-regions $\mathcal{R}_n$}{n-regions Rn}}
\begin{figure}[ht]
\centering
\begin{tikzpicture}
    \def\r{1}
    \def\e{.2cm}
    \def\a{120} 
    \def\b{90} 
    \foreach\i in {1,2,3}
        \coordinate (c\i) at (\i*\a-\a+\b:\r);
    \draw[fill=ggray, draw=white, line width=\e] (c1)--(c2)--(c3)--cycle;
    \foreach\i in {1,2,3} {
        \node[inner sep=2pt, vertexLabeled] (v\i) at (c\i) {$\i$};
    }
    \draw[thick] (v1)--(v2)--(v3)--(v1);
\end{tikzpicture}
\qquad
\begin{tikzpicture}
    \def\r{1}
    \def\e{.2cm}
    \def\a{90} 
    \def\b{45}
    \foreach\i in {1,2,3,4}
        \coordinate (c\i) at (\i*\a-\a+\b:\r);
    \draw[fill=ggray, draw=white, line width=\e] (c1)--(c2)--(c3)--(c4)--cycle;
    \draw[white, line width=\e] (c1)--(c3);
    \draw[white, line width=\e] (c2)--(c4);
    \foreach\i in {1,2,3,4} {
        \node[inner sep=2pt, vertexLabeled] (v\i) at (c\i) {$\i$};
    }
    \draw[thick] (v1)--(v2)--(v3)--(v4)--(v1);
    \draw[thick] (v1)--(v3);
    \draw[thick] (v2)--(v4);
\end{tikzpicture}
\qquad
\begin{tikzpicture}
    \def\r{1}
    \def\e{.2cm}
    \def\a{72} 
    \def\b{18}
    \foreach\i in {1,2,3,4,5}
        \coordinate (c\i) at (\i*\a-\a+\b:\r);
    \draw[fill=ggray, draw=white, line width=\e] (c1)--(c2)--(c3)--(c4)--(c5)--cycle;
    \draw[fill=ggray, draw=white, line width=\e] (c1)--(c3)--(c5)--(c2)--(c4)--cycle;
    \foreach\i in {1,2,3,4,5} {
        \node[inner sep=2pt, vertexLabeled] (v\i) at (c\i) {$\i$};
    }
    \draw[thick] (v1)--(v2)--(v3)--(v4)--(v5)--(v1);
    \draw[thick] (v1)--(v3)--(v5)--(v2)--(v4)--(v1);
\end{tikzpicture}
    \caption{The possible $n$-regions for small $n$. The shaded areas correspond to the possible choices of $a_{n+1}$ and so $r_{1,3}=1, r_{1,4}=4, r_{1,5}=11$.}
\label{fig:regionExamples}
\end{figure}
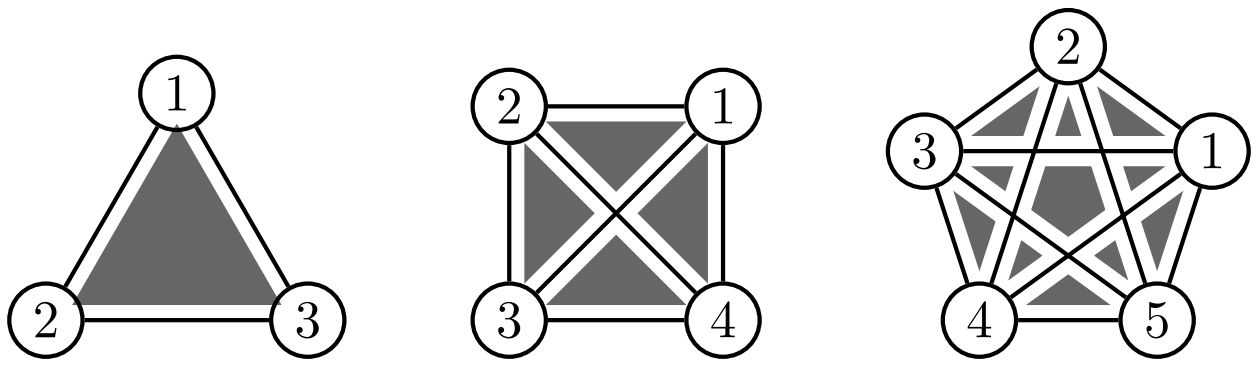

Two $n$-star configurations $A$ and $A'$ are called combinatorially equivalent, denoted by $A\sim A'$ if $\chi_{A}=\chi_{A'}$. We define the \emph{possible $n$-regions} $\mathcal{R}_n=\{[A]:A\ $n$\text{-star}\}$ as the set of all equivalence classes under $\sim$. Moreover, let $r_{1,n}=|\mathcal{R}_n|$ be their number. The first few cases are illustrated in Figure~\ref{fig:regionExamples}.

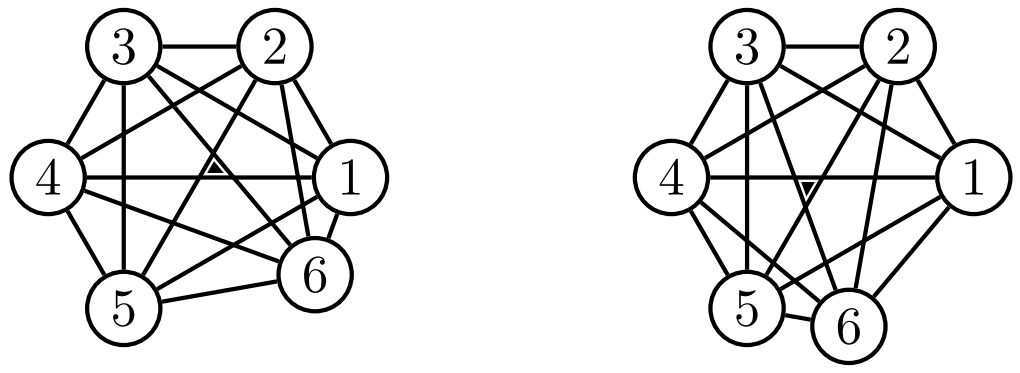
\begin{figure}[ht]
\centering
\begin{tikzpicture}
    \def\r{1}
    \def\e{.06cm}
    \def\a{60} 
    \def\b{0}
    \fill[white] (0,0) circle (\r+.3);
    \foreach\i in {1,2,3,4,5}
        \coordinate (c\i) at (\i*\a-\a+\b:\r);
    \coordinate (c6) at (6*\a-\a+\b+20:\r);
    \fill (0,0) circle(\r/2);
    \fill[white] (c1)--(c4)--(c5)--(c6)--cycle;
    \fill[white] (c1)--(c2)--(c3)--(c6)--cycle;
    \fill[white] (c2)--(c3)--(c4)--(c5)--cycle;
    \draw[draw=white, line width=\e] (c1)--(c2)--(c3)--(c4)--(c5)--(c6)--cycle;
    \draw[draw=white, line width=\e] (c1)--(c3)--(c5)--cycle;
    \draw[draw=white, line width=\e] (c2)--(c4)--(c6)--cycle;
    \draw[draw=white, line width=\e] (c1)--(c4);
    \draw[draw=white, line width=\e] (c2)--(c5);
    \draw[draw=white, line width=\e] (c3)--(c6);
    \foreach\i in {1,2,3,4,5,6}
        \node[inner sep=2pt,vertexLabeled] (v\i) at (c\i) {$\i$};
    \draw[thick] (v6)--(v1)--(v2)--(v3)--(v4)--(v5)--(v6);
    \draw[thick] (v6)--(v2)--(v4)--(v6);
    \draw[thick] (v1)--(v3)--(v5)--(v1);
    \draw[thick] (v1)--(v4);
    \draw[thick] (v2)--(v5);
    \draw[thick] (v3)--(v6);
\end{tikzpicture}
\qquad\qquad
\begin{tikzpicture}
    \def\r{1}
    \def\e{.06cm}
    \def\a{60} 
    \def\b{0}
    \fill[white] (0,0) circle (\r+.3);
    \foreach\i in {1,2,3,4,5}
        \coordinate (c\i) at (\i*\a-\a+\b:\r);
    \coordinate (c6) at (6*\a-\a+\b-20:\r);
    \fill (c1)--(c2)--(c3)--(c6)--cycle;
    \fill[white] (c1)--(c2)--(c3)--(c4)--cycle;
    \fill[white] (c1)--(c2)--(c5)--(c6)--cycle;
    \draw[draw=white, line width=\e] (c1)--(c2)--(c3)--(c4)--(c5)--(c6)--cycle;
    \draw[draw=white, line width=\e] (c1)--(c3)--(c5)--cycle;
    \draw[draw=white, line width=\e] (c2)--(c4)--(c6)--cycle;
    \draw[draw=white, line width=\e] (c1)--(c4);
    \draw[draw=white, line width=\e] (c2)--(c5);
    \draw[draw=white, line width=\e] (c3)--(c6);
    \foreach\i in {1,2,3,4,5,6}
        \node[inner sep=2pt, vertexLabeled] (v\i) at (c\i) {$\i$};
    \draw[thick] (v6)--(v1)--(v2)--(v3)--(v4)--(v5)--(v6);
    \draw[thick] (v6)--(v2)--(v4)--(v6);
    \draw[thick] (v1)--(v3)--(v5)--(v1);
    \draw[thick] (v1)--(v4);
    \draw[thick] (v2)--(v5);
    \draw[thick] (v3)--(v6);
\end{tikzpicture}
    \caption{Two different $6$-regions (shaded) that cannot coexist with the same $a_1,\ldots,a_6$. These are the only two $6$-regions that cannot coexist, making a total of $r_{1,6}=26$ regions.}
\label{fig:hexagon}
\end{figure}

Notice that for $n\geq6$, any choice of points $a_1, \ldots, a_n$ will miss out on realizing at least one possible $n$-region. See Figure~\ref{fig:hexagon}. This makes computing $r_{1,n}$ a different enumeration problem than the Moser's circle problem.

For the $n$-star configurations, the chirotope information $\chi_{A}(i<j<k)$, with $i,j,k\leq n$ is trivial, in the sense that it is always the same in all the different $n$-star configurations. This means that we are only interested in the values $\chi(i<j<n+1)$.

So, for every $n$-region $R=[A]\in \mathcal{R}_n$ we organize the relevant information of the chirotope $\chi_{A}$ in the associated skew-symmetric matrix
\[M_R=\big(\chi_{A}(i,j,n+1)\big)_{i,j\in[n]}.\]
which always satisfies the following properties:

\begin{lemma}\label{lem:simpleArgument} Let $M_R$ be the matrix associated with an $n$-region $R\in \mathcal{R}_n$. Then, for all $i,j\in[n]$ the following properties hold:
\begin{enumerate}[(i)]
    \item $(M_R)_{ij}\neq 0$ for every $i\neq j$.
    \item $(M_R)_{1n}=-1$ and for every $i\in[n-1]$, $(M_R)_{i(i+1)}=+1$.
    \item If $i<j$, $(M_R)_{ij}=+1$ and $i<k<j$, then $(M_R)_{ik},(M_R)_{kj}=+1$.
    \item If $i<j$, $(M_R)_{ij}=-1$ and ($k>j$ or $k<i$), then $(M_R)_{ik},(M_R)_{kj}=-1$.
    \item If $i>j$, $(M_R)_{ij}=-1$ and $j<k<i$, then $(M_R)_{ik},(M_R)_{kj}=-1$.
    \item If $i>j$, $(M_R)_{ij}=+1$ and ($k>i$ or $k<j$), then $(M_R)_{ik},(M_R)_{kj}=+1$.
\end{enumerate}
\end{lemma}
\begin{proof} (i) follows directly from the generality assumption and (ii) follows from the fact that $a_{n+1}$ lies inside $\conv(a_1\cdots a_n)$.

Properties (iii)-(vi) follow from one another, so it suffices to prove (iii). If $i<j$ and $(M_R)_{ij}=+1$, this means that $a_{n+1}$ is inside the convex polygon $a_1\cdots a_ia_j\cdots a_n$, which is contained in the polygon $a_1\cdots a_ia_ka_j\cdots a_n$, which proves that $(M_R)_{ik}=+1$ and $(M_R)_{kj}=+1$.
\end{proof}

Based on this remark, we transform every region $R\in\mathcal{R}_n$ injectively into a Young diagram.

\subsection{\texorpdfstring{$n$-admissible diagrams $\mathcal{Y}_n$}{n-admissible diagrams Yn}}
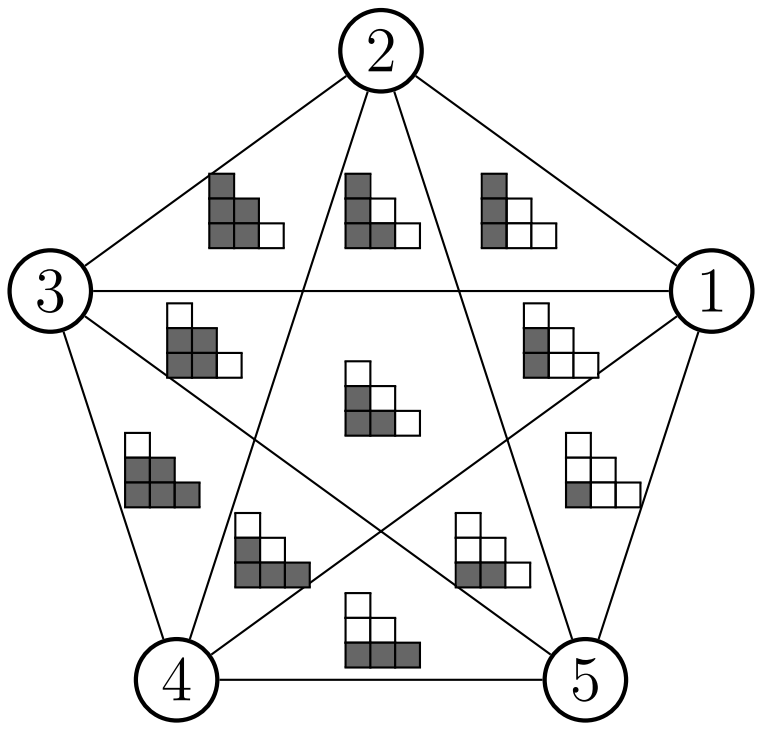
\begin{figure}
\centering
\begin{tikzpicture}
    \def\boxS{.15cm}
    \def\r{2.3}
    \def\a{72} 
    \def\b{18}
    \foreach\i in {1,2,3,4,5}
        \node[inner sep= 2pt, vertexLabeled] (v\i) at (\i*\a-\a+\b:\r) {\Large$\i$};
    \draw (v5)--(v1)--(v2)--(v3)--(v4)--(v5);
    \draw (v5)--(v2)--(v4)--(v1)--(v3)--(v5);
    \def\ri{0.53934466291} 
    \def\ro{0.66666666666} 
    \node at (18+0*72-36:\ro*\r) {\myFrenchYoung[boxSize=\boxS]{5}{1}};
    \node at (18+0*72:\ri*\r) {\myFrenchYoung[boxSize=\boxS]{5}{1,1}};
    \node at (18+1*72-36:\ro*\r) {\myFrenchYoung[boxSize=\boxS]{5}{1,1,1}};
    \node at (18+1*72:\ri*\r) {\myFrenchYoung[boxSize=\boxS]{5}{2,1,1}};
    \node at (18+2*72-36:\ro*\r) {\myFrenchYoung[boxSize=\boxS]{5}{2,2,1}};
    \node at (18+2*72:\ri*\r) {\myFrenchYoung[boxSize=\boxS]{5}{2,2}};
    \node at (18+3*72-36:\ro*\r) {\myFrenchYoung[boxSize=\boxS]{5}{3,2}};
    \node at (18+3*72:\ri*\r) {\myFrenchYoung[boxSize=\boxS]{5}{3,1}};
    \node at (18+4*72-36:\ro*\r) {\myFrenchYoung[boxSize=\boxS]{5}{3}};
    \node at (18+4*72:\ri*\r) {\myFrenchYoung[boxSize=\boxS]{5}{2}};
    \node at (0,0) {\myFrenchYoung[boxSize=\boxS]{5}{2,1}};
\end{tikzpicture}
    \caption{All the $5$-admissible diagrams (shaded boxes). The unshaded boxes stand for the entries under the diagonal that could be equal to $+1$, so the second diagonal is also removed from this notation. Notice that the last inequality $\lambda_1+k\leq 5$ excludes the two diagrams $\myFrenchYoung[boxSize=.15cm]{5}{3,2,1}$$\myFrenchYoung[boxSize=.15cm]{5}{3,1,1}$.}
\label{fig:Young}
\end{figure}

\begin{definition} We call a Young diagram associated with the partition $\lambda=(\lambda_1\geq\lambda_2\geq\cdots\geq\lambda_k>0)$ an $n$-\emph{admissible diagram}, if $\lambda$ satisfies the following properties:
\[k\geq1,\qquad \lambda_1+k\leq n,\qquad\lambda_k\leq n-k-1\]
Let $\mathcal{Y}_n$ be the set of all $n$-admissible diagrams.
\end{definition}
In particular, every $n$-admissible diagram $\lambda$ satisfies $\lambda_i\leq n-i-1$, i.e. $\lambda\leq(n-2,n-3,\ldots,1)$.
\begin{proposition}
For every $n$-region $R$, the $+1$ entries in $M_R$ below the diagonal form a Young diagram (in French notation) that is always a $n$-admissible diagram. See Figure~\ref{fig:Young} for an example.
\end{proposition}
\begin{proof} First, in order to see that the $+1$ entries in the lower left corner of the matrix $M_R$ form a Young diagram, we use Lemma~\ref{lem:simpleArgument}(vi), which in particular states that for every $+1$ under the diagonal, every entry to the left and below of this $+1$ is also equal to $+1$. Moreover, $k\geq1$ and $\lambda_k\leq n-k-1$ follow directly from Lemma~\ref{lem:simpleArgument}(ii).

For the last inequality, note that by the definition of $\lambda_1$ we have $(M_R)_{nj}=+1$ for every $j\in[\lambda_1]$. Using again Lemma~\ref{lem:simpleArgument}(vi), we get $(M_R)_{1j}=+1$ for every $j\in\{2,\ldots,\lambda_1\}$. In particular, there are at least $\lambda_1-1$ positive entries in the first row. Moreover, by the definition of $k$, we have $(M_R)_{j1}=+1$ for every $j\in\{n-k+1,\ldots,n\}$, which means that $(M_R)_{1j}=-1$ for these indices. In particular, there are exactly $k$ negative entries in the first row. So, the non-zero entries are $n-1\geq\lambda_1-1+k$, which gives the desired inequality.
\end{proof}

We finally encode each such diagram in a $0/1$ sequence based on the differences in $\lambda$.
\subsection{\texorpdfstring{$n$-admissible sequences $\mathcal{B}_n$}{n-admissible sequences Bn}}
\begin{definition}
    For $n\geq3$ we call a $0/1$ sequence $(b_{\ell})_{\ell\in[2n]}\in\{0,1\}^{2n}$ an $n$-\emph{admissible sequence} if it satisfies the following properties:
    \[b_1=1,\qquad \forall \ell\in[n]:b_{n+\ell}=1-b_{\ell},\qquad\forall \ell\in[n]:\sum_{i=1}^nb_{\ell+i}\geq1\]
    Let $\mathcal{B}_n$ be the set of all $n$-admissible sequences.
\end{definition}
Thus, an $n$-admissible sequence consists of $n$ zeros and $n$ ones, starts with $1$, is $n$-antiperiodic, and does not contain $n$ consecutive zeros.

\begin{proposition}\label{prop:YoungBraceletsBijection}
    For every $n$-admissible diagram $\lambda=(\lambda_1,\ldots,\lambda_k)\in\mathcal{Y}_n$, the following sequence is an $n$-admissible sequence:
    \[b_{\lambda}=1^{\lambda_k}01^{\lambda_{k-1}-\lambda_k}0\cdots1^{\lambda_1-\lambda_2}01^{n-k-\lambda_1}0^{\lambda_k}10^{\lambda_{k-1}-\lambda_k}1\cdots0^{\lambda_1-\lambda_2}10^{n-k-\lambda_1}\]
 Moreover, the function $\phi:\mathcal{Y}_n\to\mathcal{B}_n$ with $\phi(\lambda)=b_{\lambda}$ is a bijection.
\end{proposition}
\begin{proof}
    First, observe that all exponents are nonnegative, since $\lambda$ is a Young diagram. Moreover, $b_{\lambda}$ is well defined, since it has exactly $2n$ entries:
    \[2\big(\lambda_k+1+(\lambda_{k-1}-\lambda_k)+1+\cdots+(\lambda_1-\lambda_2)+1+(n-k-\lambda_1)\big)=2n\]
    Moreover, since $k\geq1$, $b_{\lambda}$ starts with $1$. Also, the sequence is $n$-antiperiodic by definition. Finally, for all the zeros of $b_{\lambda}$ to form a contiguous subsequence, it would be necessary that
    \[\lambda_{k-1}-\lambda_k=\cdots=\lambda_1-\lambda_2=n-k-\lambda_1=0.\]
    This would give us $\lambda_k=\lambda_{k-1}=\cdots=\lambda_1=n-k$, which is a contradiction to the definition of an $n$-admissible diagram, since $\lambda_k\leq n-k-1$.

    For the inverse, let $b\in\mathcal{B}_n$ be an $n$-admissible sequence. We define $\phi^{-1}(b)$ as follows. First we write $b$ in the form
    \[b=1^{a_k}01^{a_{k-1}}0\cdots1^{a_1}01^{a_0}0^{a_k}10^{a_{k-1}}1\cdots0^{a_1}10^{a_0}\]
    where $a_k\geq1$, $a_0,\ldots,a_{k-1}\geq0$ and $a_k+1+\cdots+a_1+1+a_0=n$ and let $\psi(b)=(a_1+\cdots+a_k,a_2+\cdots+a_k,\ldots,a_k)$ which is a Young diagram, with $k\geq1$. Also, $\lambda_1=\sum_{i=1}^ka_i=n-k-a_0\leq n-k$ and $\lambda_k=a_k\leq\sum_{i=0}^ka_k=n-k$, which can only be an equality if $a_0=a_2=\cdots=a_{k-1}=0$, which is not the case, because $b$ does not have $n$ consecutive zeros. It follows directly from the definition that $\psi\circ\phi(\lambda)=\lambda$, proving the assertion.
\end{proof}

\begin{remark}\label{rem:sequenceForm}
    It follows from the preceding proof that, for every $\lambda\in\mathcal{Y}_n$, the sequence $b_{\lambda}$ has the form
    \[b_{\lambda}=1s01^{n-k-\lambda_1}0s^*10^{n-k-\lambda_1}\]
    where $s^*$ is obtained from $s$ by interchanging zeros and ones. Moreover the part $1s0$ has exactly $k$ zeros and $\lambda_1$ ones, arranged in the form
    \[1s0=1^{\lambda_k}01^{\lambda_{k-1}-\lambda_k}0\cdots1^{\lambda_1-\lambda_2}0\]
\end{remark}

\subsection{Bijection}\label{ssec:bijection}
So far, we have established an injection $\mathcal{R}_n\to\mathcal{Y}_n$ and a bijection $\mathcal{Y}_n\leftrightarrow\mathcal{B}_n$. The aim of this section is to construct the inverse map $\mathcal{B}_n\to\mathcal{R}_n$.

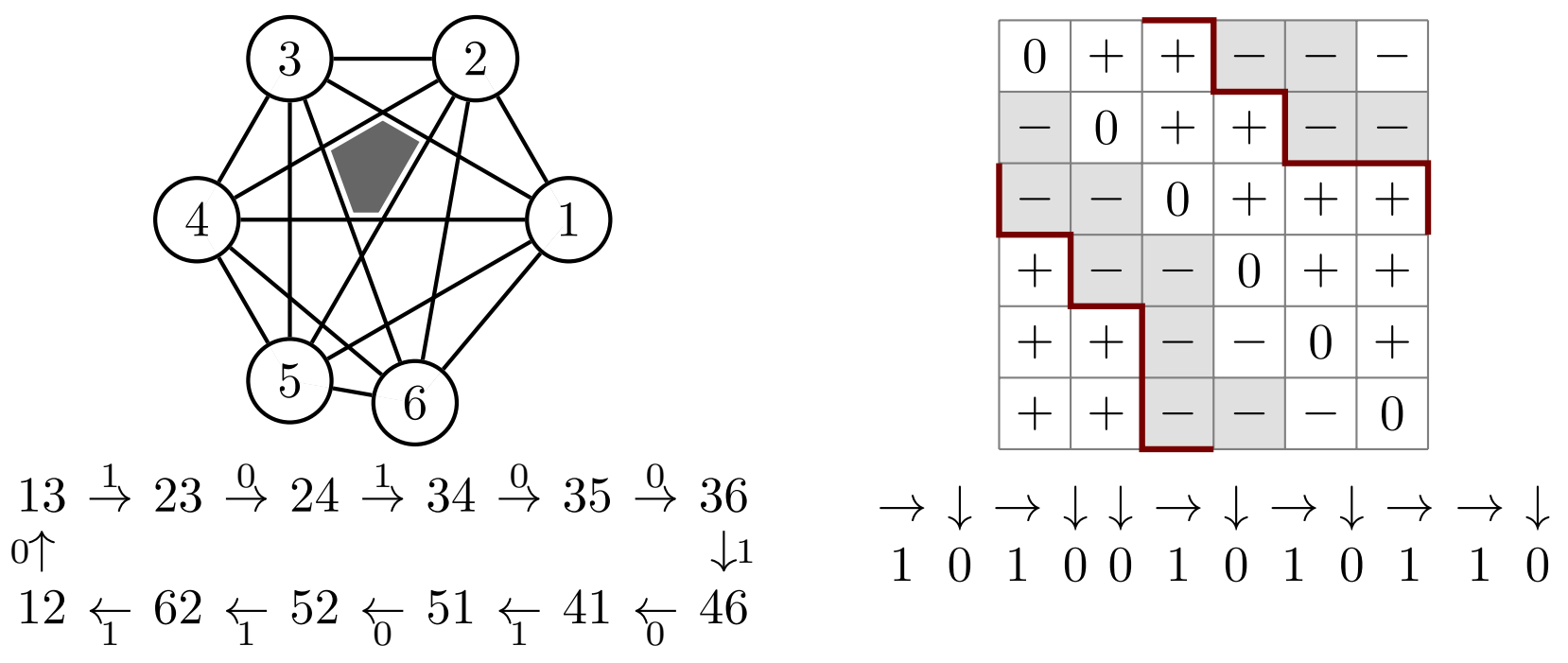
\begin{figure}[ht]
\centering
\begin{subfigure}[b]{0.45\textwidth}
    \centering
    \begin{tikzpicture}
    \def\r{1.3}
    \def\e{.1cm}
    \def\a{60} 
    \def\b{0}
    \fill[white] (0,0) circle (\r+.3);
    \foreach\i in {1,2,3,4,5}
        \coordinate (c\i) at (\i*\a-\a+\b:\r);
    \coordinate (c6) at (6*\a-\a+\b-20:\r);
    \fill (0,0) circle(\r/2);
    \fill[white] (c1)--(c4)--(c5)--(c6)--cycle;
    \fill[white] (c1)--(c2)--(c3)--(c6)--cycle;
    \fill[white] (c2)--(c3)--(c4)--(c5)--cycle;
    \draw[draw=white, line width=\e] (c1)--(c2)--(c3)--(c4)--(c5)--(c6)--cycle;
    \draw[draw=white, line width=\e] (c1)--(c3)--(c5)--cycle;
    \draw[draw=white, line width=\e] (c2)--(c4)--(c6)--cycle;
    \draw[draw=white, line width=\e] (c1)--(c4);
    \draw[draw=white, line width=\e] (c2)--(c5);
    \draw[draw=white, line width=\e] (c3)--(c6);
    \fill[ggray] (c1)--(c2)--(c3)--(c4)--(c5)--(c6)--cycle;
        \fill[white] (c1)--(c4)--(c5)--(c6)--cycle;
        \draw[draw=white, line width=\e] (c1)--(c4);
        \fill[white] (c2)--(c5)--(c6)--(c1)--cycle;
        \draw[draw=white, line width=\e] (c2)--(c5);
        \fill[white] (c3)--(c1)--(c2)--cycle;
        \draw[draw=white, line width=\e] (c3)--(c1);
        \fill[white] (c4)--(c2)--(c3)--cycle;
        \draw[draw=white, line width=\e] (c4)--(c2);
        \fill[white] (c5)--(c3)--(c4)--cycle;
        \draw[draw=white, line width=\e] (c5)--(c3);
        \fill[white] (c6)--(c3)--(c4)--(c5)--cycle;
        \draw[draw=white, line width=\e] (c6)--(c3);
    \foreach\i in {1,2,3,4,5,6}
        \node[inner sep=3pt,vertexLabeled] (v\i) at (c\i) {$\i$};
    \draw[thick] (v6)--(v1)--(v2)--(v3)--(v4)--(v5)--(v6);
    \draw[thick] (v6)--(v2)--(v4)--(v6);
    \draw[thick] (v1)--(v3)--(v5)--(v1);
    \draw[thick] (v1)--(v4);
    \draw[thick] (v2)--(v5);
    \draw[thick] (v3)--(v6);
    \end{tikzpicture}
    \begin{tikzcd}[column sep=.3cm, row sep=.3cm]
        13\ar[r,"1"]&23\ar[r,"0"]&24\ar[r,"1"]&34\ar[r,"0"]&35\ar[r,"0"]&36\ar[d,"1"]\\
        12\ar[u,"0"]&62\ar[l,"1"]&52\ar[l,"1"]&51\ar[l,"0"]&41\ar[l,"1"]&46\ar[l,"0"]
    \end{tikzcd}
    \end{subfigure}
    \begin{subfigure}[b]{0.45\textwidth}
    \centering
    \begin{tikzpicture}[scale=.5]
    \fill[ggray!20!white] (0,5)--(1,5)--(1,4)--(2,4)--(2,3)--(3,3)--(3,1)--(4,1)--(4,0)--(2,0)--(2,2)--(1,2)--(1,3)--(0,3)--cycle;
    \fill[ggray!20!white] (3,6)--(5,6)--(5,5)--(6,5)--(6,4)--(4,4)--(4,5)--(3,5)--cycle;
    \draw[step=1.0,gray,thin] (0,0) grid (6,6);
    \draw[very thick,rred] (2,6)--(3,6)--(3,5)--(4,5)--(4,4)--(5,4)--(6,4)--(6,3);
    \draw[very thick,rred] (0,4)--(0,3)--(1,3)--(1,2)--(2,2)--(2,1)--(2,0)--(3,0);
    \node at (0.5,0.5) {$+$};
    \node at (1.5,0.5) {$+$};
    \node at (2.5,0.5) {$-$};
    \node at (3.5,0.5) {$-$};
    \node at (4.5,0.5) {$-$};
    \node at (5.5,0.5) {$0$};
    \node at (0.5,1.5) {$+$};
    \node at (1.5,1.5) {$+$};
    \node at (2.5,1.5) {$-$};
    \node at (3.5,1.5) {$-$};
    \node at (4.5,1.5) {$0$};
    \node at (5.5,1.5) {$+$};
    \node at (0.5,2.5) {$+$};
    \node at (1.5,2.5) {$-$};
    \node at (2.5,2.5) {$-$};
    \node at (3.5,2.5) {$0$};
    \node at (4.5,2.5) {$+$};
    \node at (5.5,2.5) {$+$};
    \node at (0.5,3.5) {$-$};
    \node at (1.5,3.5) {$-$};
    \node at (2.5,3.5) {$0$};
    \node at (3.5,3.5) {$+$};
    \node at (4.5,3.5) {$+$};
    \node at (5.5,3.5) {$+$};
    \node at (0.5,4.5) {$-$};
    \node at (1.5,4.5) {$0$};
    \node at (2.5,4.5) {$+$};
    \node at (3.5,4.5) {$+$};
    \node at (4.5,4.5) {$-$};
    \node at (5.5,4.5) {$-$};
    \node at (0.5,5.5) {$0$};
    \node at (1.5,5.5) {$+$};
    \node at (2.5,5.5) {$+$};
    \node at (3.5,5.5) {$-$};
    \node at (4.5,5.5) {$-$};
    \node at (5.5,5.5) {$-$};
    \end{tikzpicture}
    \[\setlength{\arraycolsep}{2pt}\begin{array}{cccccccccccc}
    \rightarrow& \downarrow& \rightarrow& \downarrow& \downarrow& \rightarrow& \downarrow& \rightarrow& \downarrow& \rightarrow& \rightarrow& \downarrow\\
    1&0&1&0&0&1&0&1&0&1&1&0
    \end{array}\]
    \end{subfigure}
    \caption{Some region $R\in\mathcal{R}_6$ and its corresponding matrix $M_R$, which is encoded by the $6$-admissible diagram $\lambda_R=(2,2,1)\in\mathcal{Y}_6$. Both descriptions give us the $6$-admissible sequence $\alpha_6(R)=101001010110$. The shaded boxes $ji$ in $M_R$ correspond to the edges $ij$ produced by the algorithm in Remark~\ref{rem:algorithm}}
    \label{fig:matrixToBracelet}
\end{figure}

We first describe the map $\alpha_n:\mathcal{R}_n\to\mathcal{B}_n$  as follows: For an $n$-region $R\in\mathcal{R}_n$, $\phi(\lambda_R)$ is the path between the cells of $M_R$, which starts in the first column with a rightward step, remains between the $+1$ and $-1$ entries, and wraps around from the last row to the first. A right step is encoded by $1$ and a down step by $0$. See the right image of Figure~\ref{fig:matrixToBracelet}.

\begin{remark}\label{rem:algorithm}
The $n$-admissible sequence $\alpha_n([A])$ can also be constructed directly from the configuration $A$. We start with the edge $1j*$, where
\[j^*=\max_{j\in[n]}\{\chi_{A}(1,j,n+1)=+1\}.\]
We know that the set is non-empty, since $2$ has this property. We then rotate the edge iteratively, keeping the point $a_{n+1}$ on its left-hand side, by increasing the second index cyclically whenever possible, and otherwise increasing the first index cyclically, until we reach $1j^*$ again. So, the edge $ij$ becomes $i(j+_n1)$, if $\chi_{A}(i,j+_n1,n+1)=+1$ and $(i+_n1)j$, else, where $+_n$ stands for the sum modulo $n$ (and obviously the index $0$ stands for index $n$). Lemma~\ref{lem:simpleArgument} then guarantees that $\chi_{A}(i+_n1,j,n+1)=+1$.

If we write a $1$ every time we increase the first index and a $0$ every time we increase the second index, we retrieve $\alpha_n([A])$; see for example the left picture in Figure~\ref{fig:matrixToBracelet}.
\end{remark}
\begin{proof} By construction, whenever the algorithm increases the first index, the edge $ij$ before the increase is maximal for the first index, i.e. every number $m$ in the sequence $i+_n1,i+_n2,\ldots,i+_nn-1$ that comes before $j$ satisfies $\chi(i,m,n+1)=+1$ and every number $m$ after $j$ satisfies $\chi(i,m,n+1)=-1$. Consequently, the algorithm produces all these maximal edges for every index and thus eventually creates $1j^*$ again, which was maximal by choice. Thus, each index is increased exactly $n$ times, which gives us indeed a sequence in $\{0,1\}^{2n}$. Notice now that by the choice of the next edge, the lower left corners of the boxes $ji$ corresponding to the edges $ij$ trace exactly the corners of the path in $M_{[A]}$ that defines $b_{\lambda}$, giving us the same sequence of movements.
\end{proof}
This algorithm can easily get reversed in order to retrieve the region we started with, if we already know the positions $a_1,\ldots,a_n$, but in order to construct $\alpha_n^{-1}$ we also need to choose the positions $a_1,\ldots,a_n$ carefully, because, as we have seen, not every region exists for every choice of points $a_1,\ldots,a_n$. In fact, it is known that $n$ points in general position on the circle give rise to $\Theta(n^4)$ regions, while as we will show $r_{1,n}=\Theta(2^n)$, so as $n$ grows larger, a particular embedding of $[n]$ will only cover a minuscule fraction of all possible $n$-regions. We first need to discuss a very particular region that only exists for odd $n$.

\begin{lemma}\label{lem:centralRegion}
    Let $n=2m+1\geq3$, $a_i=(\cos\frac{2\pi i}{n},\sin\frac{2\pi i}{n})$ for all $i\in[n]$ and $a_{n+1}=(0,0)$. Then $A_n^{\circ}=\{a_1,\ldots,a_n\}$ is an $n$-star configuration. We name its region $R_n^\circ=[A_n^{\circ}]$, for which we get $\lambda_{R_n^\circ}=(m,m-1,\ldots,1)$ or equivalently $\alpha_n(R_n^{\circ})=(10)^n$. See Figure~\ref{fig:centralRegion} for an example.
\end{lemma}
\begin{proof}
    It is straightforward to verify that $A$ is an $n$-star configuration, since an odd regular polygon has no central symmetry; in particular, its center does not lie on any of its diagonals. For $\alpha_n(R_n^{\circ})$ notice that the extremal edge for every $i$ is always $i(i+_nm)$, which means that if $ij$ is an extremal edge, the next extremal edge is $(i+_n1)(j+_n1)$, i.e. we always have to increase the first position once and then the second position again once, which gives us $\alpha_n(R_n^{\circ})=(10)^n$.
\end{proof}

\begin{figure}[ht]
    \centering
\begin{subfigure}[b]{.45\textwidth}
    \centering
    \begin{tikzpicture}[baseline]
    \def\r{1.3}
    \def\e{.1cm}
    \def\a{51.4285714286} 
    \foreach\i in {1,...,7}
        \coordinate (c\i) at (\i*\a-\a-12.85:\r);
    \fill[ggray] (c1)--(c2)--(c3)--(c4)--(c5)--(c6)--(c7)--cycle;
    \draw[draw=white, line width=\e] (c1)--(c4)--(c7)--(c3)--(c6)--(c2)--(c5)--cycle;
    \fill[white] (c1)--(c2)--(c3)--(c4)--cycle;
    \fill[white] (c2)--(c3)--(c4)--(c5)--cycle;
    \fill[white] (c3)--(c4)--(c5)--(c6)--cycle;
    \fill[white] (c4)--(c5)--(c6)--(c7)--cycle;
    \fill[white] (c5)--(c6)--(c7)--(c1)--cycle;
    \fill[white] (c6)--(c7)--(c1)--(c2)--cycle;
    \fill[white] (c7)--(c1)--(c2)--(c3)--cycle;
    \foreach\i in {1,...,7}
        \node[vertexLabeled, inner sep=3pt] (v\i) at (c\i) {\small$\i$};
    \draw[thick] (v1)--(v2)--(v3)--(v4)--(v5)--(v6)--(v7)--(v1);
    \draw[thick] (v1)--(v4)--(v7)--(v3)--(v6)--(v2)--(v5)--(v1);
\end{tikzpicture}
\[R_7^{\circ}\]
\end{subfigure}
\begin{subfigure}[b]{.45\textwidth}
    \[\myFrenchYoung[boxSize=.45cm]{7}{3,2,1}\]
    \vspace*{.5em}
    \[\lambda_{R_7^{\circ}}=(3,2,1)\qquad\ \]
\end{subfigure}
\caption{The case $\alpha_7(R_7^{\circ})=10101010101010$.}
\label{fig:centralRegion}
\end{figure}

\begin{proposition}\label{prop:bijection}
    The map $\alpha_n:\mathcal{R}_n\to\mathcal{B}_n$ defined above is a bijection.
\end{proposition}
\begin{proof}
    We now construct the inverse of $\alpha_n$ for an arbitrary $b\in\mathcal{B}_n$. We first express $b$ in block form as follows:
    \[1^{a_1}0^{b_1}1^{a_2}0^{a_3}\cdots0^{b_p}1^{a_{p+1}}\]
    where $a_1,\ldots,a_p,b_1,\ldots,b_p\geq1$ and $a_{p+1}\geq0$. Let $\ell\in[p]$ be the index of the block of zeros containing $b_{n+1}=0$. Then, because of the $n$-antiperiodicity of $b$ it is true that $a_{p+1}+a_1=b_{\ell}$ and $a_i=b_{\ell+_p(i-1)}$ for every $i\in[p]\setminus\{1\}$. Also, because of the $n$-antiperiodicity, $b_p=a_{\ell}=b_{2\ell-1}$, i.e. $p=2\ell-1$ is an odd number. Also, since $b$ does not have a consecutive subsequence of $n$ zeros, $p\neq1$.
    We label the $n$ occurrences of $1$ consecutively, from $1$ to $n$ and define the function $c:[n]\to[p]$ assigning to each occurrence of $1$ the index of the block to which it belongs, with the last $a_{p+1}$ ones getting assigned to the first block.
    
    Let then $b^{\circ}=(10)^p\in\mathcal{B}_p$. We have established in Lemma~\ref{lem:centralRegion} that there exists a $p$-star configuration $A_p^{\circ}$ such that $\alpha_p([A_p^{\circ}])=b^{\circ}$. Let $\{w_1,\ldots,w_{p+1}\}$ be this configuration. We now define a new configuration $A'=\{v_1,v_2,\ldots,v_n,v_{n+1}\}$ such that $v_{n+1}=w_{p+1}$ and $v_1,\ldots,v_n$ are in convex position, labeled counterclockwise, in such a way that for every $i\in[n]$ $\|v_i-w_{c(i)}\|<\varepsilon$ for some small $\varepsilon>0$. For example, suppose that $w_1,\ldots,w_p$ are the points on the circle constructed in Lemma~\ref{lem:centralRegion}, we can for every $j\in[p]$ distribute the $|c^{-1}(j)|$ corresponding points around $w_j$ on the circle, such that they remain on the arc centered at $w_j$ of the length $2\varepsilon$. See Figure~\ref{fig:constructionInverse} for an example of this construction.
    
    The parameter $\varepsilon$ should be chosen sufficiently small so that the movement of the edge $v_iv_j$ to $w_{c(i)}w_{c(j)}$ does not cross the point $v_{n+1}=w_{p+1}$, or equivalently $\chi_{A'}(i,j,n+1)=\chi_{A_p^{\circ}}(c(i),c(j),p+1)$.

    It remains to show that $\alpha([A'])=b$. If we fix a $t\in[p]$, notice that all extremal edges starting at some $i$ with $c(i)=t$ have the same endpoint, namely $j^*=\max c^{-1}(t+_p\ell)$, where the maximum is on the cyclical interval. This means that if $ij^*$ is the extremal edge of $i$, we rotate it by increasing its first coordinate until this first coordinate has a different value under $c$. It follows that
    \[\alpha_n([A'])=1^{a_1}0^{b_1'}1^{a_2}0^{b_2'}\cdots0^{b_p'}1^{a_{p+1}}\]
    for some $b_1',b_2',\ldots,b_p'\geq1$, which means that $\alpha_n([A'])$ and $b$ agree on the sizes of the ones blocks. Because both are $n$-antiperiodic this also fixes the block sizes of the zeros, thus we have $\alpha_n([A'])=b$.
\end{proof}

\begin{figure}[ht]
    \centering
\begin{subfigure}[b]{.45\textwidth}
    \centering
\begin{tikzpicture}
    \def\r{1.3}
    \def\e{.1cm}
    \def\a{51.4285714286} 
    \foreach\i in {1,...,7}
        \coordinate (c\i) at (\i*\a-\a-12.85:\r);
    \fill[ggray] (c1)--(c2)--(c3)--(c4)--(c5)--(c6)--(c7)--cycle;
    \draw[draw=white, line width=\e] (c1)--(c4)--(c7)--(c3)--(c6)--(c2)--(c5)--cycle;
    \fill[white] (c1)--(c2)--(c3)--(c4)--cycle;
    \fill[white] (c2)--(c3)--(c4)--(c5)--cycle;
    \fill[white] (c3)--(c4)--(c5)--(c6)--cycle;
    \fill[white] (c4)--(c5)--(c6)--(c7)--cycle;
    \fill[white] (c5)--(c6)--(c7)--(c1)--cycle;
    \fill[white] (c6)--(c7)--(c1)--(c2)--cycle;
    \fill[white] (c7)--(c1)--(c2)--(c3)--cycle;
    \foreach\i in {1,...,7}
        \node[vertexLabeled, inner sep=3pt] (v\i) at (c\i) {\small$\i$};
    \draw[thick] (v1)--(v2)--(v3)--(v4)--(v5)--(v6)--(v7)--(v1);
    \draw[thick] (v1)--(v4)--(v7)--(v3)--(v6)--(v2)--(v5)--(v1);
\end{tikzpicture}
\[\setlength{\arraycolsep}{0.3pt}\begin{array}{ccccccccccccccc}b^{\circ}=\,&1&0&1&0&1&0&1&0&1&0&1&0&1&0\\&1&&2&&3&&4&&5&&6&&7\end{array}\]
\end{subfigure}
\begin{subfigure}[b]{.45\linewidth}
\centering
\begin{tikzpicture}
    \def\r{1.5}
    \def\R{1.7}
    \def\a{51.4285714286} 
    \def\b{-12.85} 
    \def\eps{6}
    \def\bigEps{13}
    \coordinate (c12) at (0*\a+\b-\eps:\r);
    \node (n12) at (0*\a+\b-\bigEps:\R) {\small$12$};
    \coordinate (c1) at (0*\a+\b+0:\r);
    \node at (0*\a+\b+0:\R) {\small$1$};
    \coordinate (c2) at (0*\a+\b+\eps:\r);
    \node at (0*\a+\b+\bigEps:\R) {\small$2$};
    \coordinate (c3) at (1*\a+\b-\eps/2:\r);
    \node at (1*\a+\b-\bigEps/2:\R) {\small$3$};
    \coordinate (c4) at (1*\a+\b+\eps/2:\r);
    \node at (1*\a+\b+\bigEps/2:\R) {\small$4$};
    \coordinate (c5) at (2*\a+\b:\r);
    \node at (2*\a+\b:\R) {\small$5$};
    \coordinate (c6) at (3*\a+\b-\eps:\r);
    \node at (3*\a+\b-\bigEps:\R) {\small$6$};
    \coordinate (c7) at (3*\a+\b:\r);
    \node at (3*\a+\b:\R) {\small$7$};
    \coordinate (c8) at (3*\a+\b+\eps:\r);
    \node at (3*\a+\b+\bigEps:\R) {\small$8$};
    \coordinate (c9) at (4*\a+\b:\r);
    \node at (4*\a+\b:\R) {\small$9$};
    \coordinate (c10) at (5*\a+\b:\r);
    \node at (5*\a+\b:\R) {\small$10$};
    \coordinate (c11) at (6*\a+\b:\r);
    \node at (6*\a+\b:\R) {\small$11$};
    \fill[ggray] (c1)--(c2)--(c3)--(c4)--(c5)--(c6)--(c7)--(c8)--(c9)--(c10)--(c11)--(c12)--cycle;
    \fill[white] (c1)--(c9)--(c10)--(c11)--(c12)--cycle;
    \fill[white] (c2)--(c9)--(c10)--(c11)--(c12)--(c1)--cycle;
    \fill[white] (c3)--(c10)--(c11)--(c12)--(c1)--(c2)--cycle;
    \fill[white] (c4)--(c10)--(c11)--(c12)--(c1)--(c2)--(c3)--cycle;
    \fill[white] (c5)--(c11)--(c12)--(c1)--(c2)--(c3)--(c4)--cycle;
    \fill[white] (c6)--(c12)--(c1)--(c2)--(c3)--(c4)--(c5)--cycle;
    \fill[white] (c7)--(c12)--(c1)--(c2)--(c3)--(c4)--(c5)--(c6)--cycle;
    \fill[white] (c8)--(c12)--(c1)--(c2)--(c3)--(c4)--(c5)--(c6)--(c7)--cycle;
    \fill[white] (c9)--(c3)--(c4)--(c5)--(c6)--(c7)--(c8)--cycle;
    \fill[white] (c10)--(c5)--(c6)--(c7)--(c8)--(c9)--cycle;
    \fill[white] (c11)--(c6)--(c7)--(c8)--(c9)--(c10)--cycle;
    \fill[white] (c12)--(c9)--(c10)--(c11)--cycle;
    \draw[white, line width=.1cm] (c1)--(c9);
    \draw[white, line width=.1cm] (c2)--(c9);
    \draw[white, line width=.1cm] (c3)--(c10);
    \draw[white, line width=.1cm] (c4)--(c10);
    \draw[white, line width=.1cm] (c5)--(c11);
    \draw[white, line width=.1cm] (c6)--(c12);
    \draw[white, line width=.1cm] (c7)--(c12);
    \draw[white, line width=.1cm] (c8)--(c12);
    \draw[white, line width=.1cm] (c9)--(c3);
    \draw[white, line width=.1cm] (c10)--(c5);
    \draw[white, line width=.1cm] (c11)--(c6);
    \draw[white, line width=.1cm] (c12)--(c9);
    \draw[thick] (c1)--(c2)--(c3)--(c4)--(c5)--(c6)--(c7)--(c8)--(c9)--(c10)--(c11)--(c12)--cycle;
    \draw[thick] (c1)--(c9);
    \draw[thick] (c2)--(c9);
    \draw[thick] (c2)--(c10);
    \draw[thick] (c3)--(c10);
    \draw[thick] (c4)--(c10);
    \draw[thick] (c4)--(c11);
    \draw[thick] (c5)--(c11);
    \draw[thick] (c5)--(c12);
    \draw[thick] (c6)--(c12);
    \draw[thick] (c7)--(c12);
    \draw[thick] (c8)--(c12);
    \draw[thick] (c8)--(c1);
    \draw[thick] (c8)--(c2);
    \draw[thick] (c8)--(c3);
    \draw[thick] (c9)--(c3);
    \draw[thick] (c9)--(c4);
    \draw[thick] (c9)--(c5);
    \draw[thick] (c10)--(c5);
    \draw[thick] (c10)--(c6);
    \draw[thick] (c11)--(c6);
    \draw[thick] (c11)--(c7);
    \draw[thick] (c11)--(c8);
    \draw[thick] (c11)--(c9);
    \draw[thick] (c12)--(c9);
\end{tikzpicture}
\[\setlength{\arraycolsep}{0.3pt}\begin{array}{ccccccccccccccccccccccccc}b=
&1&1&0&1&1&0&1&0&1&1&1&0&0&0&1&0&0&1&0&1&0&0&0&1\\
&1&1& &2&2& &3& &4&4&4& & & &5& & &6& &7& & & &1\end{array}\]
\end{subfigure}
\caption{Creating one $12$-star configuration from a particular $b\in\mathcal{B}_{12}$. For the construction we use the $R_7^{\circ}$, since $b$ has $7$ blocks of ones. Below $b$ we write for every $i\in[12]$ the value $c(i)$ corresponding to the $i$-th occurrence of $1$ inside $b$.}
\label{fig:constructionInverse}
\end{figure}

\begin{theorem}\label{thm:regions}
    Let $n\geq3$. Then
    \begin{equation}\label{eq:numRegions}
        r_{1,n}=2^{n-1}-n.
    \end{equation}
\end{theorem}
\begin{proof}
 Due to the bijection proven in Proposition~\ref{prop:bijection} we know that $r_{1,n}=|\mathcal{R}_n|=|\mathcal{B}_n|=2^{n-1}-n$. Indeed, every $n$-admissible sequence starts with an $1$, the next $n-1$ entries can be chosen freely and the rest is then fixed because of the $n$-antiperiodicity. We must nonetheless subtract from $2^{n-1}$ the $n$ sequences that have $n$ consecutive zeros.
\end{proof}
Sequence~\ref{eq:numRegions} is also known as the Eulerian numbers $A(n-1,1)$, see \href{https://oeis.org/A000295}{[A000295]} for a list of other possible combinatorial interpretations.

\subsection{Action of the Dihedral group}\label{ssec:action}
The dihedral group
\[D_n=\left<r,f:r^n=f^2=(fr)^2=e\right>\]
acts naturally on $\mathcal{R}_n$.
\begin{definition} 
Let $n\geq3$. We define $\rho_n:D_n\times\mathcal{R}_n\to\mathcal{R}_n$ to be the following action. For $R=[A]\in\mathcal{R}_n$ let $r\cdot[A]=[B]$, with $b_i=a_{i+_n1}$ and $f\cdot[A]=[C]$, with $c_i=a_{n-i}$. Equivalently, $M_{r\cdot R}$ arises from $M_R$ by rotating its first column and row to be the last column and row. Moreover $M_{f\cdot R}$ arises from $M_R$ by reversing the order of its rows and columns (i.e. reflecting across the antidiagonal of the matrix).
\end{definition}
Note that the reflection of $\lambda\in\mathcal{Y}_n$ corresponds to taking its conjugate partition $\lambda^T$.

\begin{lemma}\label{lem:action}
    Under the bijection $\alpha_n:\mathcal{R}_n\to\mathcal{B}_n$, $\mathcal{B}_n$ inherits the following $D_n$-action. Let $b\in\mathcal{B}_n$. For the rotation $r$, write $b$ in the form $10^a1t$, where $a\geq0$ and $t$ is the end of the sequence. Then $r\cdot b=1t10^a$. For the flip, write $b$ in the form $s10^b$, where $b\geq0$ and $s$ is the start of the sequence. Then $f\cdot b=1\bar{s}0^b$, where $\bar{s}$ is $s$ reversed.
\end{lemma}
\begin{proof} Notice that the path described in Subsection~\ref{ssec:bijection} in $M_{r\cdot R}$ is obtained from the path encoded by $b$ in $M_R$ by moving the first $1$, together with the zeros (if any exist) preceding the next $1$, to the end of the sequence.
For the flipped matrix $M_{f\cdot R}$ the path gets reversed and rotated cyclically. Thus, it remains to identify which occurrence of $1$ in $b$ becomes the first entry of $f\cdot b$. Notice that inside $M_{f\cdot R}$ the first $\rightarrow$ move corresponds to the last $\downarrow$ move in $M_R$, before the path cyclically moves from the $n$-th row to the first, so the last $\downarrow$ move among its first $n$ moves. This corresponds to the last $0$ among $b_1,b_2,\ldots,b_n$, whose corresponding $1$ is $n$ places afterwards. This means that after $b$ gets reversed we have to bring this particular $1$ to the first place by rotating the preceding zeros (if any exist) to the end.
\end{proof}

One can argue that the number of orbits of this action is the true number of combinatorially different $n$-star configurations, which can be attained by the Burnside's Lemma. For enumeration of such combinatorial objects we refer the reader to \cite{bracelets}, where one can also find the proof of Proposition~\ref{prop:orbits}.
\begin{definition}
    Let $n\geq3$. Then an $n$-admissible bracelet is a black-white bracelet with $2n$ beads that is \emph{antisymmetric} (also called \emph{self-dual}) and contains no monochromatic $n$-bead chain. Let $\tilde{\mathcal{B}}_n$ be the set of all $n$-admissible bracelets (both rotating and flipping is allowed). See Figure~\ref{fig:bracelets} for a list of every element in $\tilde{\mathcal{B}}_n$ for the first few $n$.
\end{definition}
Every element of $\mathcal{B}_n$ is mapped to a bracelet in $\tilde{\mathcal{B}}_n$ by joining the two ends of the sequence and replacing each $1$ with a black bead and each $0$ with a white bead. Notice that from Lemma~\ref{lem:action} follows that $|\mathcal{R}_n/D_n|=|\tilde{\mathcal{B}}_n|$.
\begin{proposition}[\cite{bracelets}]\label{prop:orbits}
    Let $n\geq3$. Then the number of combinatorially different $n$-star configurations under the action of $D_n$  equals to:
        \[|\mathcal{R}_n/D_n|=|\tilde{\mathcal{B}}_n|=\frac{1}{2}\left(2^{\left\lfloor\frac{n-1}{2}\right\rfloor}+\sum_{2\nmid d\mid n}\frac{\phi(d)2^{\frac{n}{d}}}{2n}\right)-1\]
(see \href{https://oeis.org/A007147}{[A007147]}).
\end{proposition}

\subsection{\texorpdfstring{The $n$ excluded configurations}{The n excluded configurations}}
\begin{figure}[ht]
\centering
\begin{tikzpicture}
    \def\boxS{.15cm}
    \def\r{2.3}
    \def\a{72} 
    \def\b{18}
    \foreach\i in {1,2,3,4,5}
        \node[inner sep= 2pt, vertexLabeled] (v\i) at (\i*\a-\a+\b:\r) {\Large$\i$};
    \draw (v5)--(v1)--(v2)--(v3)--(v4)--(v5);
    \draw (v5)--(v2)--(v4)--(v1)--(v3)--(v5);
    \def\ri{0.53934466291} 
    \def\ro{0.66666666666} 
    \node at (18+0*72-36:\ro*\r) {\myFrenchYoungNoStair[boxSize=\boxS]{5}{1}};
    \node at (18+0*72-36:\r) {$\emptyset$};
    \node at (18+0*72:\ri*\r) {\myFrenchYoungNoStair[boxSize=\boxS]{5}{1,1}};
    \node at (18+1*72-36:\ro*\r) {\myFrenchYoungNoStair[boxSize=\boxS]{5}{1,1,1}};
    \node at (18+1*72-36:\r) {\myFrenchYoungNoStair[boxSize=\boxS]{5}{1,1,1,1}};
    \node at (18+1*72:\ri*\r) {\myFrenchYoungNoStair[boxSize=\boxS]{5}{2,1,1}};
    \node at (18+2*72-36:\ro*\r) {\myFrenchYoungNoStair[boxSize=\boxS]{5}{2,2,1}};
    \node at (18+2*72-36:\r) {\myFrenchYoungNoStair[boxSize=\boxS]{5}{2,2,2}};
    \node at (18+2*72:\ri*\r) {\myFrenchYoungNoStair[boxSize=\boxS]{5}{2,2}};
    \node at (18+3*72-36:\ro*\r) {\myFrenchYoungNoStair[boxSize=\boxS]{5}{3,2}};
    \node at (18+3*72-36:\r) {\myFrenchYoungNoStair[boxSize=\boxS]{5}{3,3}};
    \node at (18+3*72:\ri*\r) {\myFrenchYoungNoStair[boxSize=\boxS]{5}{3,1}};
    \node at (18+4*72-36:\ro*\r) {\myFrenchYoungNoStair[boxSize=\boxS]{5}{3}};
    \node at (18+4*72-36:\r) {\myFrenchYoungNoStair[boxSize=\boxS]{5}{4}};
    \node at (18+4*72:\ri*\r) {\myFrenchYoungNoStair[boxSize=\boxS]{5}{2}};
    \node at (0,0) {\myFrenchYoungNoStair[boxSize=\boxS]{5}{2,1}};
\end{tikzpicture}
    \caption{The same as Figure~\ref{fig:Young}, but this time the diagrams only need satisfy $\lambda_1+k\leq5$.}
\label{fig:YoungPlusPlus}
\end{figure}
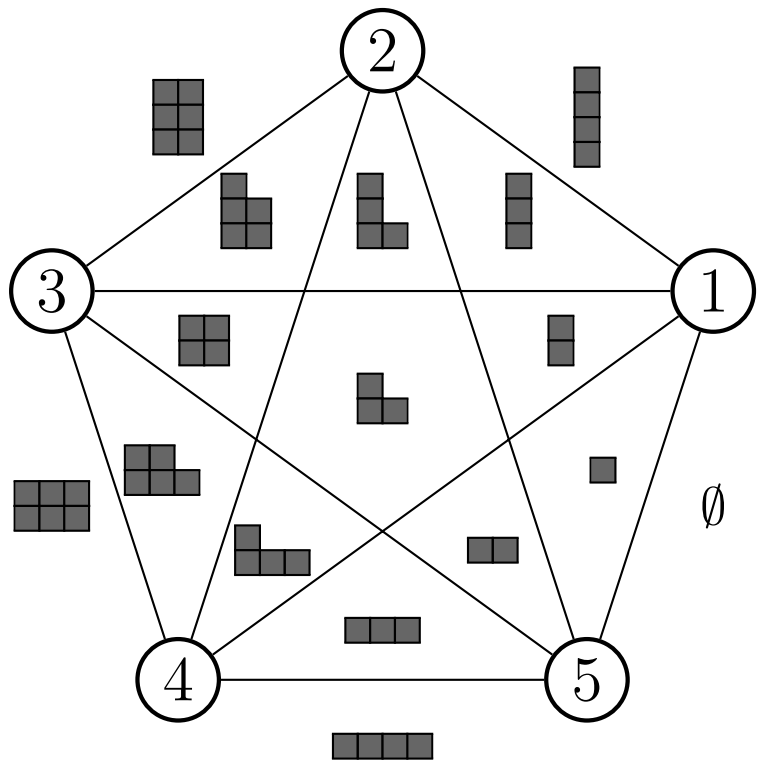

Before concluding this section, we discuss the $n$ elements excluded from $\mathcal{B}_n$ and show that they have a corresponding combinatorial meaning in terms of both $\mathcal{R}_n$ and $\mathcal{Y}_n$. In Proposition~\ref{prop:YoungBraceletsBijection}, while proving the bijection between $\mathcal{Y}_n$ and $\mathcal{B}_n$, we saw that for an $n$-sequence to have $n$ consecutive zeros is equivalent either with a diagram of the form $\lambda_1=\lambda_2=\cdots=\lambda_k$ with $\lambda_1+k=n$, or with $k=0$ (the empty diagram). There are exactly $n$ such diagrams, one for each value of $k=0,\ldots,n-1$. It turns out that these $n$ diagrams correspond to the $n$ combinatorial configurations, in which $a_{n+1}$ is not an interior point but instead lies just outside the edge $i(i+_n1)$:
\begin{proposition}
    Let $A=(a_1,a_2,\ldots,a_n,a_{n+1})$ be a configuration in which all points are in a convex position, the points with the labels $1,\ldots,n$ lie in counterclockwise order, and the last point lies cyclically between $k\in[n-1]$ and $k+1$. Then:
    \[\chi_{A}(i<j<n+1)=\left\{\begin{array}{ll}-1&,i<=k<j\\+1&else\end{array}\right.\]
    If $a_{n+1}$ lies between $a_n$ and $a_1$, then $\chi_A(i<j<n+1)=+1$ for every pair $i,j\in[n]$. The Young diagrams inside $M_R$ are then exactly the $n$ rectangular diagrams with $\lambda_1+k=n$.
\end{proposition}
\begin{proof}
 Since all points are in convex position $\chi(\{i_1,i_2,i_3\})=+1$ if we trace the indices $i_1,i_2,i_3$ counterclockwise. So, for the case of $i_2=n+1$, using the antisymmetry of $\chi_A$ we get $\chi_A(i<j<n+1)=-1$, whilst for $i_1=n+1$ or $i_3=n+1$ the value of $\chi_A$, after permuting, remains $+1$. If $a_{n+1}$ lies between $a_n$ and $a_1$, the sign is trivially positive for every pair, since all points are labeled in counterclockwise order.
\end{proof}
The set containing $\mathcal{Y}_n$ and these additional diagrams consists of Young diagrams $(\lambda_1,\ldots,\lambda_k)$ that only satisfy $\lambda_1+k\leq n$.
See Figure~\ref{fig:YoungPlusPlus} for an example of these additional regions and diagrams.

\begin{remark}
    In this light, Theorem~\ref{thm:regions} can be written as
    \[r_{1,n}=2^{n-1}-nr_{0,n+1}\]
    where $r_{0,n+1}=1$ is the number of the different combinatorial regions with $n+1$ points, none of them are in the interior and the coefficient $n$ in front of it counts in how many ways the last point can be part of the boundary. 
\end{remark}

\section{Enumerating Triangulations}
Throughout the paper, we use the notation
\[C_n=\frac{1}{n+1}\binom{2n}{n}\]
for the $n$-th Catalan numbers.
\subsection{Enumeration for a fixed region}
\begin{definition}
For $n\geq3$ and $\lambda\in\mathcal{Y}_n$, let $t_{n,\lambda}=|\mathcal{T}^{\text{full}}_A|$ denote the number of full triangulations of an $n$-star configuration $A$ satisfying $\lambda_{[A]}=\lambda$.
\end{definition}
In \cite{wheels}, the authors provide a method for counting the full triangulations of an $n$-star configuration.
They show that, when the interior point lies sufficiently close to an edge, the number of full triangulations equals $C_{n-1}-C_{n-2}$ (Lemma~1).
As the point moves through the interior, crossing a chord changes the number of full triangulations. Orienting the chord so that the point crosses it from left to right, the corresponding increase is given by $C_{p-1}C_{n-p}-C_{p}C_{n-p-1}$ (Lemma~2), where $p+1$ is the number of vertices of the left polygon that this chord creates by splitting the bigger polygon. The right polygon has then $(n-p+1)$ vertices.
We have adjusted the indices in both lemmas to match our notation.

In terms of Young diagrams, crossing an edge from left to right corresponds to adding a box to the diagram. Thus, the Young diagram provides a combinatorial encoding of a sequence of chord crossings encountered as the interior point moves through the interior of the $n$-gon. This allows us to translate the geometric counting argument of \cite{wheels} into a formula depending only on the associated partition.
\begin{lemma}\label{lem:inductionBox}
Let $n\geq3$ and let $\lambda,\lambda'\in\mathcal{Y}_n$ differ by exactly one box in row $s\in[k]$. That is, $\lambda_s'=\lambda_s+1$ and for every other index $i\neq s$, $\lambda_i'=\lambda_i$. Then
\[t_{n,\lambda'}=t_{n,\lambda}+\big(C_{(\lambda_s+s)-1}C_{n-(\lambda_s+s)}-C_{\lambda_s+s}C_{n-(\lambda_s+s)-1}\big)\]
\end{lemma}
\begin{proof}
    By the construction of $\lambda$, the point $a_{n+1}$ lies to the left of the edge $a_{n-s+1}a_{\lambda_s}$ and on the right side of the edge $a_{n-s+1}a_{\lambda_s+1}$. Similarly, $a_{n+1}'$ lies to the left of $a_{n-s+1}a_{\lambda_s+1}$ and on the right side of $a_{n-s+1}a_{\lambda_s+2}$. Thus, in this transition, the interior point crosses the edge $a_{\lambda_s+1}a_{n-s+1}$ from the left to the right. The polygon to the left of this edge has vertices 
    \[a_{n-s+1},a_{n-s},\ldots,a_n,a_1,\ldots,a_{\lambda_s+1},\]
    which are $s+\lambda_s+1=p+1$ many. This gives $p=\lambda_s+s$, which proves the lemma.
\end{proof}
Observe that adding several boxes to the same row gives rise to a telescoping sum. This observation allows us to express the number of full triangulations of an $n$-star configuration entirely in terms of its associated Young diagram.
\begin{theorem}\label{thm:enumTriangFromYoung}
    Let $n\geq3$ and $\lambda\in\mathcal{Y}_n$. Then
    \begin{equation}\label{eq:enumTriangFromYoung}
        t_{n,\lambda}=\sum_{i=1}^k\Big(C_{i-1}C_{n-i}-C_{(\lambda_i+i)-1}C_{n-(\lambda_i+i)}\Big)
    \end{equation}
\end{theorem}
\begin{proof}
    We prove the result by induction on the Young poset, starting from the minimal element $\lambda_{\min}=(1)$. Notice that if $\lambda\in\mathcal{Y}_n$ then every element of the interval $[(1),\lambda]$ inside the Young poset is also an element of $\mathcal{Y}_n$. So, the minimum $\lambda_{\min}$ corresponds to the inner point being ``very close'' to the edge $n1$. By the first lemma of \cite{wheels}, $t_{n,(1)}=C_{n-1}-C_{n-2}$, which agrees with Equation~\ref{eq:enumTriangFromYoung}. Assume that the formula holds for some $\lambda$. We use Lemma~\ref{lem:inductionBox} to determine the number of triangulations associated with $\lambda'$, obtained by adding one box to $\lambda$ in row $s$. That is $\lambda_s'=\lambda_s+1$ and otherwise $\lambda_i'=\lambda_i$.
    \begin{align*}
        t_{n,\lambda'}&=t_{n,\lambda}+\big(C_{(\lambda_s+s)-1}C_{n-(\lambda_s+s)}-C_{\lambda_s+s}C_{n-(\lambda_s+s)-1}\big)\\
        &=\sum_{i=1}^k\Big(C_{i-1}C_{n-i}-C_{(\lambda_i+i)-1}C_{n-(\lambda_i+i)}\Big)+\\
        &\qquad\qquad\qquad\qquad+\big(C_{(\lambda_s+s)-1}C_{n-(\lambda_s+s)}-C_{\lambda_s+s}C_{n-(\lambda_s+s)-1}\big)\\
        &=\sum_{i=1}^k\Big(C_{i-1}C_{n-i}-C_{\lambda_i'+i-1}C_{n-\lambda_i'-i}\Big)\qedhere
    \end{align*}
\end{proof}

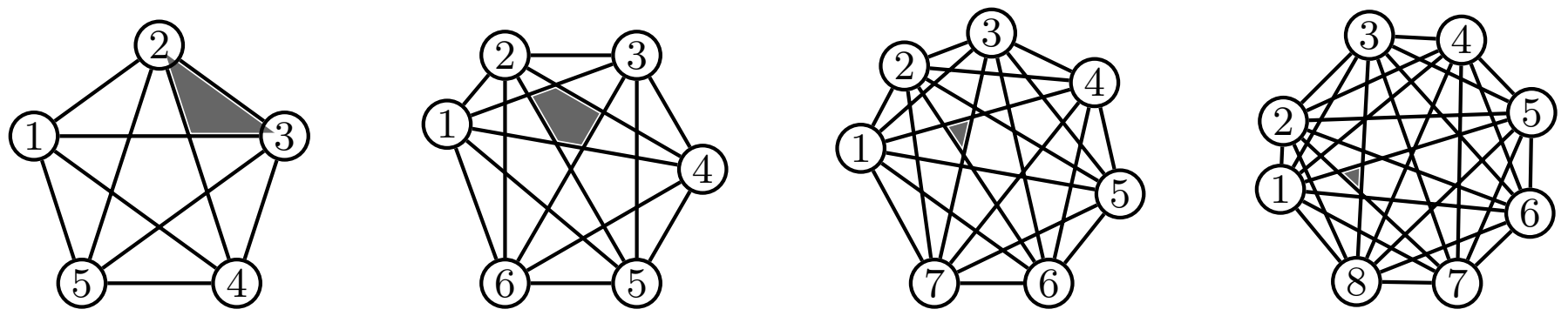
\begin{figure}[ht]
    \centering
    \begin{subfigure}[b]{.24\textwidth}
    \centering
    \begin{tikzpicture}
        \def\b{0.17cm}
        \def\r{1}
        \def\e{.05cm}
        \def\a{72} 
        \foreach\i in {1,2,3,4,5}
            \coordinate (c\i) at (90-\i*\a+\a+\a:\r);
        \fill[ggray] (c1)--(c2)--(c3)--(c4)--(c5)--cycle;
        \draw[draw=white, line width=\e] (c1)--(c2)--(c3)--(c4)--(c5)--cycle;
        \draw[draw=white, line width=\e] (c1)--(c3)--(c5)--(c2)--(c4)--cycle;
        \fill[white] (c3)--(c4)--(c5)--(c1)--cycle;
        \fill[white] (c2)--(c4)--(c5)--(c1)--cycle;
        \foreach\i in {1,2,3,4,5}
            \node[vertexLabeled, inner sep=1pt, minimum size=5pt] (v\i) at (c\i) {\small$\i$};
        \draw[thick] (v5)--(v1)--(v2)--(v3)--(v4)--(v5);
        \draw[thick] (v5)--(v2)--(v4)--(v1)--(v3)--(v5);
    \end{tikzpicture}
    \end{subfigure}
    \begin{subfigure}[b]{.24\textwidth}
    \centering
    \begin{tikzpicture}
        \def\r{1}
        \def\e{.05cm}
        \def\a{60} 
        \coordinate (c1) at (60-\a+\a+2*\a-20:\r);
        \foreach\i in {2,3,4,5,6}
            \coordinate (c\i) at (60-\i*\a+\a+2*\a:\r);
        \fill[ggray] (c1)--(c2)--(c3)--(c4)--(c5)--(c6)--cycle;
        \draw[draw=white, line width=\e] (c1)--(c2)--(c3)--(c4)--(c5)--(c6)--cycle;
        \draw[draw=white, line width=\e] (c1)--(c3)--(c5)--cycle;
        \draw[draw=white, line width=\e] (c2)--(c4)--(c6)--cycle;
        \draw[draw=white, line width=\e] (c1)--(c4);
        \draw[draw=white, line width=\e] (c2)--(c5);
        \draw[draw=white, line width=\e] (c3)--(c6);
        \fill[white] (c1)--(c3)--(c2)--cycle;
        \fill[white] (c1)--(c4)--(c5)--(c6)--cycle;
        \fill[white] (c2)--(c4)--(c3)--cycle;
        \fill[white] (c2)--(c5)--(c6)--(c1)--cycle;
        \fill[white] (c3)--(c6)--(c5)--(c4)--cycle;
        \foreach\i in {1,2,3,4,5,6}
            \node[vertexLabeled, inner sep=1pt, minimum size=5pt] (v\i) at (c\i) {\small$\i$};
        \draw[thick] (v6)--(v1)--(v2)--(v3)--(v4)--(v5)--(v6);
        \draw[thick] (v5)--(v1)--(v3)--(v5);
        \draw[thick] (v2)--(v4)--(v6)--(v2);
        \draw[thick] (v1)--(v4);
        \draw[thick] (v2)--(v5);
        \draw[thick] (v3)--(v6);
    \end{tikzpicture}
    \end{subfigure}
    \begin{subfigure}[b]{.24\textwidth}
    \centering
    \begin{tikzpicture}
        \def\r{1}
        \def\e{.05cm}
        \def\a{51.4285714286} 
        \coordinate (c1) at (90-\a+\a+2*\a-20:\r);
        \coordinate (c2) at (90-2*\a+\a+2*\a-10:\r);
        \foreach\i in {3,4,5,6,7}
            \coordinate (c\i) at (90-\i*\a+\a+2*\a:\r);
        \fill[ggray] (c1)--(c2)--(c3)--(c4)--(c5)--(c6)--(c7)--cycle;
        \draw[draw=white, line width=\e] (c1)--(c2)--(c3)--(c4)--(c5)--(c6)--(c7)--cycle;
        \draw[draw=white, line width=\e] (c1)--(c3)--(c5)--(c7)--(c2)--(c4)--(c6)--cycle;
        \draw[draw=white, line width=\e] (c1)--(c4)--(c7)--(c3)--(c6)--(c2)--(c5)--cycle;
        \fill[white] (c1)--(c4)--(c3)--(c2)--cycle;
        \fill[white] (c1)--(c5)--(c6)--(c7)--cycle;
        \fill[white] (c2)--(c5)--(c4)--(c3)--cycle;
        \fill[white] (c2)--(c6)--(c7)--(c1)--cycle;
        \fill[white] (c3)--(c7)--(c6)--(c5)--(c4)--cycle;
        \fill[white] (c4)--(c7)--(c6)--(c5)--cycle;
        \foreach\i in {1,2,3,4,5,6,7}
            \node[vertexLabeled, inner sep=1pt, minimum size=5pt] (v\i) at (c\i) {\small$\i$};
        \draw[thick] (v7)--(v1)--(v2)--(v3)--(v4)--(v5)--(v6)--(v7);
        \draw[thick] (v7)--(v2)--(v4)--(v6)--(v1)--(v3)--(v5)--(v7);
        \draw[thick] (v7)--(v3)--(v6)--(v2)--(v5)--(v1)--(v4)--(v7);
    \end{tikzpicture}
    \end{subfigure}
    \begin{subfigure}[b]{.24\textwidth}
    \centering
    \begin{tikzpicture}
        \def\r{1}
        \def\e{.05cm}
        \def\a{45} 
        \coordinate (c1) at (45+22.5-\a+\a+3*\a-10:\r);
        \coordinate (c2) at (45+22.5-2*\a+\a+3*\a+5:\r);
        \coordinate (c3) at (45+22.5-3*\a+\a+3*\a-5:\r);
        \foreach\i in {4,5,6,7,8}
            \coordinate (c\i) at (44+22.5-\i*\a+\a+3*\a:\r);
        \fill[ggray] (c1)--(c2)--(c3)--(c4)--(c5)--(c6)--(c7)--(c8)--cycle;
        \draw[draw=white, line width=\e] (c1)--(c2)--(c3)--(c4)--(c5)--(c6)--(c7)--(c8)--cycle;
        \draw[draw=white, line width=\e] (c1)--(c3)--(c5)--(c7)--cycle;
        \draw[draw=white, line width=\e] (c2)--(c4)--(c6)--(c8)--cycle;
        \draw[draw=white, line width=\e] (c1)--(c4)--(c7)--(c2)--(c5)--(c8)--(c3)--(c6)--cycle;
        \draw[draw=white, line width=\e] (c1)--(c5);
        \draw[draw=white, line width=\e] (c2)--(c6);
        \draw[draw=white, line width=\e] (c3)--(c7);
        \draw[draw=white, line width=\e] (c4)--(c8);
        \fill[white] (c1)--(c5)--(c4)--(c3)--(c2)--cycle;
        \fill[white] (c1)--(c6)--(c7)--(c8)--cycle;
        \fill[white] (c2)--(c6)--(c5)--(c4)--(c3)--cycle;
        \fill[white] (c2)--(c7)--(c8)--(c1)--cycle;
        \fill[white] (c3)--(c8)--(c7)--(c6)--(c5)--(c4)--cycle;
        \foreach\i in {1,2,3,4,5,6,7,8}
            \node[vertexLabeled, inner sep=1pt, minimum size=5pt] (v\i) at (c\i) {\small$\i$};
        \draw[thick] (v1)--(v2)--(v3)--(v4)--(v5)--(v6)--(v7)--(v8)--(v1);
        \draw[thick] (v1)--(v3)--(v5)--(v7)--(v1);
        \draw[thick] (v2)--(v4)--(v6)--(v8)--(v2);
        \draw[thick] (v1)--(v4)--(v7)--(v2)--(v5)--(v8)--(v3)--(v6)--(v1);
        \draw[thick] (v1)--(v5);
        \draw[thick] (v2)--(v6);
        \draw[thick] (v3)--(v7);
        \draw[thick] (v4)--(v8);
    \end{tikzpicture}
    \end{subfigure}
    \caption{The regions for increasing $n$ that correspond to the fixed diagram $\lambda=(3,2)$. The number of triangulations for these configurations is given by $t_{n,(3,2)}=C_0C_{n-1}-C_3C_{n-4}+C_1C_{n-2}-C_3C_{n-4}=C_{n-1}+C_{n-2}-10C_{n-4}$ for $n\geq5$.}
    \label{fig:fixedYoung}
\end{figure}

A consequence of Theorem~\ref{thm:enumTriangFromYoung} is a similar effect that the authors of \cite{wheels} noticed for $\lambda=(1)$ holds for every non-empty Young diagram $\lambda$. Observe first that there always exists an integer  $N\in\mathbb{Z}_{\geq3}$ such that $\lambda\in\mathcal{Y}_n$ if and only if $n\geq N$. Equation~\ref{eq:enumTriangFromYoung} then becomes for every $n$ the same linear combination of (at most) $2k$ terms of the form $C_{n-i}$. See, for example, Figure~\ref{fig:fixedYoung} for such a sequence.

\subsection{Enumeration for a fixed dihedral orbit}
In Section~\ref{ssec:action}, we described the action of the dihedral group on $\mathcal{R}_n$ and established a correspondence between its orbits and the bracelets in $\tilde{\mathcal{B}}_n$. We now express the number of full triangulations associated with a dihedral orbit of an $n$-region entirely in terms of its combinatorial representation, namely, the corresponding bracelet.

\begin{definition}
    Let $n\geq3$ and $\tilde{b}\in\tilde{\mathcal{B}}_n$. Since $t_{n,\lambda}=t_{n,\lambda'}$ whenever $\lambda,\lambda'\in\mathcal{Y}_n$ correspond to the same bracelet, we denote by $t_{n,\tilde{b}}=t_{n,\lambda}$ the common number of full triangulations associated with the regions in the given orbit $\tilde{b}$.
\end{definition}

\begin{figure}[ht!]
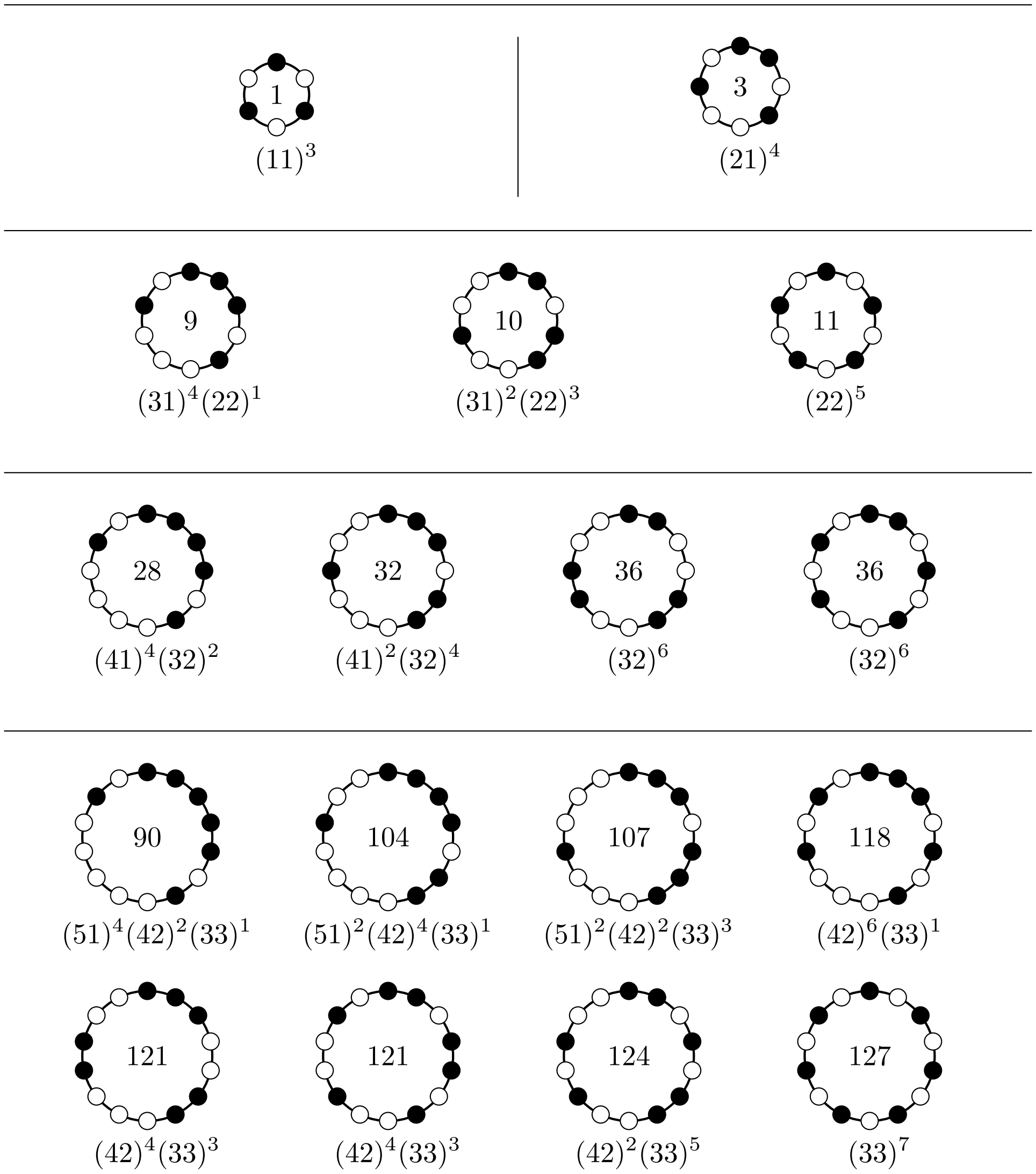

    \def\beadSize{6pt}
    \hrulefill\\[1em]
    \centering
    \begin{subfigure}[b]{.45\textwidth}
        \centering
        \myBracelet[beadSize=\beadSize, braceletRadius=.4cm, labelInside={$1$}]{1,0,1}
        \vspace*{-.5em}
        \[(11)^3\]
    \end{subfigure}\vrule depth 1em
    \begin{subfigure}[b]{.45\textwidth}
        \centering
        \myBracelet[beadSize=\beadSize, braceletRadius=.5cm, labelInside={$3$}]{1,1,0,1}
        \vspace*{-.5em}
        \[(21)^4\]
    \end{subfigure}\\[1em]\hrulefill\\[1em]
%
    \begin{subfigure}[b]{.3\textwidth}
        \centering
        \myBracelet[beadSize=\beadSize, braceletRadius=.6cm, labelInside={$9$}]{1,1,1,0,1}
        \vspace*{-.5em}
        \[(31)^4(22)^1\]
    \end{subfigure}
    \begin{subfigure}[b]{.3\textwidth}
        \centering
        \myBracelet[beadSize=\beadSize, braceletRadius=.6cm, labelInside={$10$}]{1,1,0,1,1}
        \vspace*{-.5em}
        \[(31)^2(22)^3\]
    \end{subfigure}
    \begin{subfigure}[b]{.3\textwidth}
        \centering
        \myBracelet[beadSize=\beadSize, braceletRadius=.6cm, labelInside={$11$}]{1,0,1,0,1}
        \vspace*{-.5em}
        \[(22)^5\]
    \end{subfigure}\\[1em]\hrulefill\\[1em]
    \begin{subfigure}[b]{.225\textwidth}
        \centering
        \myBracelet[beadSize=\beadSize, braceletRadius=.7cm, labelInside={$28$}]{1,1,1,1,0,1}
        \vspace*{-.5em}
        \[(41)^4(32)^2\]
    \end{subfigure}
    \begin{subfigure}[b]{.225\textwidth}
        \centering
        \myBracelet[beadSize=\beadSize, braceletRadius=.7cm, labelInside={$32$}]{1,1,1,0,1,1}
        \vspace*{-.5em}
        \[(41)^2(32)^4\]
    \end{subfigure}
    \begin{subfigure}{.225\textwidth}
        \centering
        \myBracelet[beadSize=\beadSize, braceletRadius=.7cm, labelInside={$36$}]{1,1,0,0,1,1}
        \vspace*{-.5em}
        \[(32)^6\]
    \end{subfigure}
    \begin{subfigure}{.225\textwidth}
        \centering
        \myBracelet[beadSize=\beadSize, braceletRadius=.7cm, labelInside={$36$}]{1,1,0,1,0,1}
        \vspace*{-.5em}
        \[(32)^6\]
    \end{subfigure}\\[1em]\hrulefill\\[1em]
    \begin{subfigure}{.225\textwidth}
        \centering
        \myBracelet[beadSize=\beadSize, braceletRadius=.8cm, labelInside={$90$}]{1,1,1,1,1,0,1}
        \vspace*{-.5em}
        \[(51)^4(42)^2(33)^1\]
    \end{subfigure}
    \begin{subfigure}{.225\textwidth}
        \centering
        \myBracelet[beadSize=\beadSize, braceletRadius=.8cm, labelInside={$104$}]{1,1,1,1,0,1,1}
        \vspace*{-.5em}
        \[(51)^2(42)^4(33)^1\]
    \end{subfigure}
    \begin{subfigure}{.225\textwidth}
        \centering
        \myBracelet[beadSize=\beadSize, braceletRadius=.8cm, labelInside={$107$}]{1,1,1,0,1,1,1}
        \vspace*{-.5em}
        \[(51)^2(42)^2(33)^3\]
    \end{subfigure}
    \begin{subfigure}{.225\textwidth}
        \centering
        \myBracelet[beadSize=\beadSize, braceletRadius=.8cm, labelInside={$118$}]{1,1,1,0,1,0,1}
        \vspace*{-.5em}
        \[(42)^6(33)^1\]
    \end{subfigure}\\[1em]
    \begin{subfigure}{.225\textwidth}
    \centering
        \myBracelet[beadSize=\beadSize, braceletRadius=.8cm, labelInside={$121$}]{1,1,1,0,0,1,1}
        \vspace*{-.5em}
        \[(42)^4(33)^3\]
    \end{subfigure}
    \begin{subfigure}{.225\textwidth}
    \centering
        \myBracelet[beadSize=\beadSize, braceletRadius=.8cm, labelInside={$121$}]{1,1,0,1,1,0,1}
        \vspace*{-.5em}
        \[(42)^4(33)^3\]
    \end{subfigure}
    \begin{subfigure}{.225\textwidth}
    \centering
        \myBracelet[beadSize=\beadSize, braceletRadius=.8cm, labelInside={$124$}]{1,1,0,1,0,1,1}
        \vspace*{-.5em}
        \[(42)^2(33)^5\]
    \end{subfigure}
    \begin{subfigure}{.225\textwidth}
    \centering
        \myBracelet[beadSize=\beadSize, braceletRadius=.8cm, labelInside={$127$}]{1,0,1,0,1,0,1}
        \vspace*{-.5em}
        \[(33)^7\]
    \end{subfigure}
    \caption{All bracelets for $3\leq n\leq7$. Below each bracelet, we display the multiset of its antipodal cuts; the corresponding value of $t_{n,\tilde{b}}$ is shown in the center.}
    \label{fig:bracelets}
\end{figure}

\begin{definition}
    Let $n\geq3$ and $\tilde{b}\in\tilde{\mathcal{B}}_n$ be a bracelet. In such a bracelet, each bead has an antipodal bead of the opposite color. Each antipodal pair determines $(n_1,n_2)$ an integer partition of $n-1$, which we call an \emph{antipodal cut}, where $n_1$ and $n_2$ count the black beads on the two (unordered) arcs between the beads of the pair. Let $\Pi_{\tilde{b}}$ be the multiset that contains all antipodal cuts. There are $n$ of them, since $\tilde{b}$ has exactly $n$ antipodal pairs. See Figure~\ref{fig:bracelets} for the values of $\Pi_{\tilde{b}}$ and $t_{1,\tilde{b}}$ for some small $n$.
\end{definition}

\begin{theorem}\label{thm:enumTriangFromBra} Let $n\geq3$ and $\tilde{b}\in\tilde{\mathcal{B}}_n$. Then
\begin{equation}
    t_{n,\tilde{b}}=\frac{1}{2}\left(C_n-\sum_{(n_1,n_2)\in\Pi_{\tilde{b}}}C_{n_1}C_{n_2}\right)
\end{equation}
\end{theorem}
\begin{proof}
    Choose $\lambda\in\mathcal{Y}_n$ such that $\tilde{b}=[b_{\lambda}]$. Since the conjugate partition $\lambda^T=(\lambda'_1,\ldots,\lambda'_{k'})$ corresponds to the reflected sequence $f\cdot b_{\lambda}$, it determines the same bracelet. Consequently, $t_{n,\tilde{b}}=\frac{1}{2}\left(t_{n,\lambda}+t_{n,\lambda^T}\right)$.
    We first combine the positive terms in Equation~\ref{eq:enumTriangFromYoung} for both $\lambda$ and its conjugate partition $\lambda^T$.
    Observe that for $k'$ (the number of rows of $\lambda^T$) holds $k'=\lambda_1$. Thus $\sum_{i=1}^{k'}C_{i-1}C_{n-i}=\sum_{i=1}^{\lambda_1}C_{i-1}C_{n-i}$ and by the change of variables $j=n-i$ this equals $\sum_{n-\lambda_1}^{n-1}C_jC_{n-j-1}$. Reindexing the sum corresponding to $\lambda$ as well so that it starts at $i=0$, the positive terms of the sum $t_{n,\lambda}+t_{n,\lambda^T}$ are
    \begin{align*}
    &\sum_{i=1}^kC_{i-1}C_{n-i}+\sum_{i=1}^{k'}C_{i-1}C_{n-i}=\sum_{i=0}^{k-1}C_iC_{n-1-i}+\sum_{i=n-\lambda_1}^{n-1}C_iC_{n-1-i}\\
    =&\sum_{i=0}^{n-1}C_iC_{n-1-i}-\sum_{i=k}^{n-\lambda_1-1}C_iC_{n-1-i}=C_n-\sum_{i=k}^{n-\lambda_1-1}C_iC_{n-1-i}
    \end{align*}
    Notice that since $k+\lambda_1\leq n$ we can always combine the two sums into a full Catalan convolution, with a middle range of $n-k-\lambda_1$ indices removed. Consequently, $t_{n,\tilde{b}}=\frac{1}{2}\left(C_n-N\right)$, where
    \begin{equation}\label{eq:negativeParts}
        N = \sum_{i=1}^kC_{(\lambda_i+i)-1}C_{n-(\lambda_i+i)}+\sum_{i=1}^{k'}C_{(\lambda'_i+i)-1}C_{n-(\lambda'_i+i)}+\sum_{i=k}^{n-\lambda_1-1}C_iC_{n-1-i}
    \end{equation}
    It remains to show that the three sums in the expression for $N$ enumerate all antipodal cuts of $\tilde{b}$, with the corresponding Catalan weights.
    We use Remark~\ref{rem:sequenceForm} to write $b_{\lambda}$ in the form $1s01^{n-k-\lambda_1}0s^*10^{n-k-\lambda_1}$, where $s^*$ is the dual sequence of $s$, obtained from $s$ by interchanging the zeros and the ones. The starting block $1s0$ has exactly $k$ zeros. The starting block of the second half $0s^*1$ has exactly $\lambda_1$ zeros, and the last block $0^{n-k-\lambda_1}$ has the rest $n-k-\lambda_1$ zeros. We investigate the three cases separately, where the $0$ (white bead) of an antipodal pair lies in one of these three blocks of $b_{\lambda}$.
    \begin{itemize}
        \item First, consider the $(k-i+1)$-st zero in the first block, indicated by the underline.
        \[1s0=1^{\lambda_k}01^{\lambda_{k-1}-\lambda_k}0\cdots1^{\lambda_i-\lambda_{i+1}}\underline{0}\cdots1^{\lambda_1-\lambda_2}0.\]
        Cutting the bracelet along this $0$ and its antipodal $1$ divides it into two arcs, one of which is:
        \[0^{\lambda_{i-1}-\lambda_{i-2}}1\cdots0^{\lambda_1-\lambda_2}10^{n-k-\lambda_1}1^{\lambda_k}01^{\lambda_{k-1}-\lambda_k}0\cdots1^{\lambda_i-\lambda_{i+1}}\]
        This arc contains exactly $i-1+\lambda_i$ black beads. The other arc therefore contains the rest $n-1-(i-1+\lambda_i)$ black beads. From this, it follows that the $k$ antipodal cuts in this first case is the multiset
        \[\left\{\big((\lambda_1+i)-1,n-(\lambda_i+i)\big)\in\Pi_{\tilde{b}}: i\in[k]\right\}.\]
        \item Next, considerl the $\lambda_1$ zeros in the second block $0s^*1$. By the description of the dihedral action in Section~\ref{ssec:action}, 
        \[b_{\lambda^T}=f\cdot b_{\lambda}=1\bar{s}^*01^{n-k-\lambda_1}0\bar{s}10^{n-k-\lambda_1}\]
        where $\bar{s}^*$ and $s^*$ have both exactly $\lambda_1$ zeros. This means that we can apply the same analysis as in the first case, only this time for $\lambda^T$.
        So, the $\lambda_1$ antipodal cuts of $\tilde{b}$ along a $0$ inside $0s^*1$ is the multiset
        \[\left\{\big((\lambda'_i+i)-1,n-(\lambda'_i+i)\big)\in\Pi_{\tilde{b}}:i\in[k']\right\}.\]
        \item Finally, for the last $r:=n-k-\lambda_1$ antipodal pairs, consider the $(r-i+1)$-st zero in the last block $0^r$ for some $i\in[r]$. In this case, one of the two arcs between the selected $0$ and its antipodal $1$ is:
        \[0^{i-1}1s01^{r-i}\]
        for every $i\in[n-k-\lambda_1]$. This contains exactly $\lambda_1+r-i=n-k-i$ black beads and so the other arc contains $k+i-1$ black beads. From this it follows that the last $n-k-\lambda_1$ antipodal cuts is the multiset 
        \[\left\{\big(k+i-1,n-k-i\big)\in\Pi_{\tilde{b}}:i\in[n-k-\lambda_1]\}\right\}.\]
    \end{itemize}
    The three cases exhaust all $n$ antipodal cuts in $\Pi_{\tilde{b}}$ and correspond precisely to the three sums in Equation~\ref{eq:negativeParts}. Hence,
    \[N=\sum_{(n_1,n_2)\in\Pi_{\tilde{b}}}C_{n_1}C_{n_2},\]
    finishing the proof of the theorem.
\end{proof}

\subsection{Enumerating topological triangulations}
In all our discussion so far, we have assumed the points to have a particular embedding, but a natural question is also how many different triangulations there exist if we let the points move freely in the plane under the constraint that the points labeled by $1$ up to $n$ stay in convex position cyclically ordered, the point labeled by $n+1$ stays strictly in their convex hull, and that no three points are collinear.

\begin{figure}[ht]
    \centering
    \begin{tikzpicture}
        \def\R{2.5cm}
        \fill[rred] (0,0) circle (\R);
        \draw[white, line width=.1cm] (0,0) circle (\R);
        \draw[thick] (0,0) circle (\R);
        \coordinate (c1) at (32:\R);
        \coordinate (c2) at (111:\R);
        \coordinate (c3) at (194:\R);
        \coordinate (c4) at (-104:\R);
        \coordinate (c5) at (-22:\R);
        \coordinate (cs) at (\R/4,\R/3);
        \fill[white] (c1)--(c2)--(c3)--(c4)--(c5)--cycle;
        \draw[white, line width=.1cm] (c1)--(c2)--(c3)--(c4)--(c5)--cycle;
        \node[fill=white, vertexLabeled, inner sep=1pt] (a1) at (c1) {$a_{p_1}$};
        \node[fill=white, vertexLabeled, inner sep=1pt] (a2) at (c2) {$a_{p_2}$};
        \node[fill=white, vertexLabeled, inner sep=1pt] (a3) at (c3) {$a_{p_3}$};
        \node[fill=white, vertexLabeled, inner sep=1pt] (a4) at (c4) {$a_{p_4}$};
        \node[fill=white, vertexLabeled, inner sep=1pt] (a5) at (c5) {$a_{p_5}$};
        \node[fill=white, vertexLabeled, inner sep=1pt] (as) at (cs) {$a_{n+1}$};
        \draw[thick] (a1)--(a2)--(a3)--(a4)--(a5)--(a1);
        \draw[thick] (as)--(a1);
        \draw[thick] (as)--(a2);
        \draw[thick] (as)--(a3);
        \draw[thick] (as)--(a4);
        \draw[thick] (as)--(a5);
    \end{tikzpicture}
\caption{A sketch of how every triangulation $T$ with $P_T=\{p_1<\cdots<p_5\}$ looks like. The shaded areas are convex polygons that each such $T$ needs to triangulate.}
\label{fig:polygon}
\end{figure}
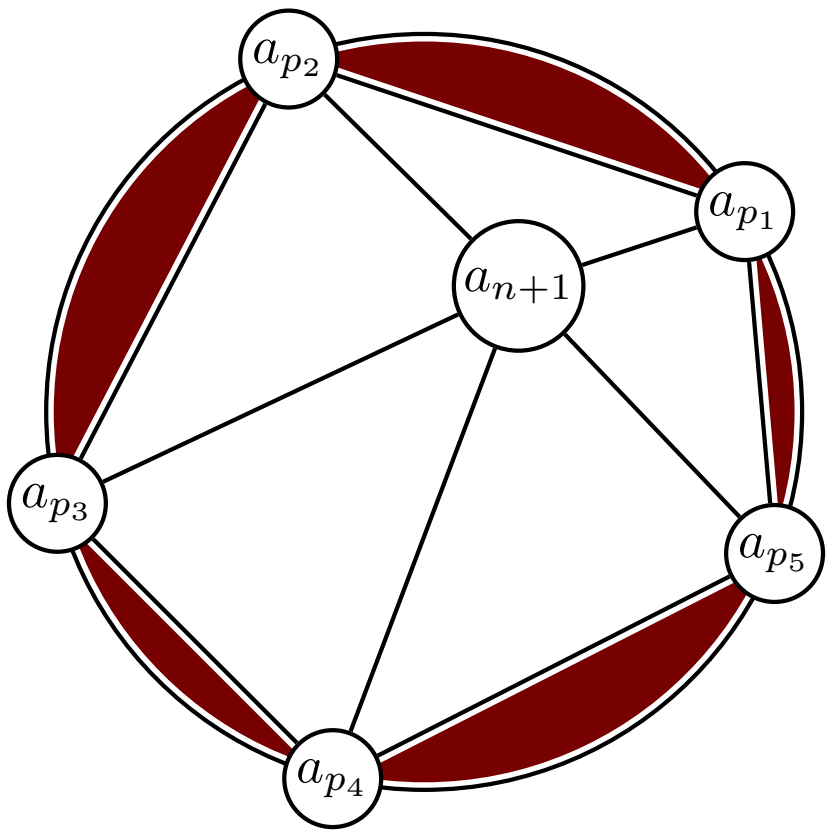

\begin{definition}
    Let $n\geq3$. Then let
    \[\mathcal{T}_{1,n}^{\text{full}}=\bigcup_{A\ n\text{-star}}\mathcal{T}^{\text{full}}_A\]
    be the set of all \emph{topological star triangulations} and $t_{1,n}=\left|\mathcal{T}_{1,n}^{\text{full}}\right|$ be its size.
\end{definition}
Notice that $t_{1,n}$ does not simply equal the sum of all $t_{n,\lambda}$, since every topological star triangulation exists for multiple $\lambda\in\mathcal{Y}_n$. We base our computation of $t_{1,n}$ on the following observation.
\begin{lemma}\label{lem:polygon}
    Let $n\geq3$, $A$ be an $n$-star configuration, and $T\in\mathcal{T}^{\text{full}}_A$. Moreover, let $P_T\subseteq[n]$ be the set of indices such that $a_{n+1}a_p$ is an edge in $T$ if and only if $p\in P$. Then, for every fixed set $P=\{p_1<\cdots<p_k\}\subseteq[n]$ with $k\geq3$, the following holds.
        \[\left|\big\{T:T\in\mathcal{T}^{\text{full}}_{1,n},\ P_T=P\big\}\right|=C_{p_2-p_1-1}\cdots C_{p_k-p_{k-1}-1}C_{p_1+n-p_k-1}\]
\end{lemma}
\begin{proof}
    Since $a_{n+1}$ is only connected to the points $\{a_{p_i}:i\in[k]\}$ within some triangulation $T$ with $P_T=\{p_1<\cdots p_k\}$, the edge $a_{p_i}a_{p_{i+_k1}}$ should also be an edge of $T$ for every $i\in[k]$. This only leaves a choice for every such $T$ on how the convex polygons $\{p_i,p_i+_n1,\ldots,p_{i+_k1}\}$ are triangulated, each of them having exactly $C_{p_{i+1}-p_i-1}$ different triangulations, if $i\in[k-1]$ and the last having $C_{p_1+n-p_k-1}$ triangulations. See also Figure~\ref{fig:polygon}
\end{proof}
Since we are considering topological star triangulations, the inner point can be in the interior of any polygon defined by boundary vertices and fixing such a polygon always produces different triangulations, since $P_{T'}\neq P_T$ means that $T'\neq T$. So, from Lemma~\ref{lem:polygon} follows the following equation.
\begin{equation}\label{eq:productCatalan}
    t_{1,n}=\sum_{k=3}^n\ \sum_{1\leq p_1<p_2<\cdots<p_k\leq n}C_{p_2-p_1-1}\cdots C_{p_k-p_{k-1}-1}C_{p_1+n-p_k-1}
\end{equation}

\begin{theorem}\label{thm:enumTriangTop}
    Let $n\geq3$. Then
    \begin{equation}\label{eq:topTr}
        t_{1,n}=\binom{2n-3}{n-3}
    \end{equation}
\end{theorem}
\begin{proof}
    Notice that choosing some $\{p_1<\cdots<p_k\}\subseteq[n]$ is equivalent to the following choices. First we choose some $q\in[n]$, then some $d_i\geq0$ for every $i\in[k-1]$ and finally some $d_k\geq p_1-1$, such that $\sum_{i=1}^kd_i=n-k$. The bijection between these choices and the indices $p_i$ is given by $q=p_1$, $d_i=p_{i+1}-p_i$ for every $i\in[k-1]$ and $d_k=p_1+n-p_k-1$. Using this change of variables in Equation~\ref{eq:productCatalan}, we get:
    \[t_{1,n}=\sum_{k=3}^n\ \sum_{\substack{d_1,\dots,d_k\geq0\\ q\in[d_k+1]\\\sum _{i=1}^kd_i=n-k}}C_{d_1}C_{d_2}\cdots C_{d_k}=\sum_{\substack{d_1,\dots,d_k\geq0\\\sum _{i=1}^kd_i=n-k}}(d_k+1)C_{d_1}\cdots C_{d_k}\]
    Let us now define the generating function $T_1(x)=\sum_{n\geq3}t_{1,n}x^n$. Using the expression above for $t_{1,n}$ gives:
    \begin{align*}
        T_1(x)&=\sum_{k\geq3}\ \sum_{d_1,\ldots,d_k\geq0}(d_k+1)C_{d_1}\cdots C_{d_k}x^{d_1+\cdots+d_k+k}\\
        &=\sum_{k\geq3}\left(\sum_{d\geq0}C_dx^{d+1}\right)^{k-1}\sum_{d_k\geq0}(d_k+1)C_{d_k}x^{d_k+1}\\
        &=\sum_{k\geq3}\big(xC(x)\big)^{k-1}\big(xC(x)+x^2C'(x)\big)\\
        &=\frac{x^3C(x)^2\big(C(x)+xC'(x)\big)}{1-xC(x)}
    \end{align*}
    where $C(x)=\frac{1-\sqrt{1-4x}}{2x}$ is the generating function of the Catalan numbers. Substituting the explicit expression of $C(x)$ and simplifying, we obtain
    \begin{equation}\label{eq:genFunctionTop}
        T_1(x)=\frac{1-3x}{2\sqrt{1-4x}}+\frac{1}{2}x-\frac{1}{2}
    \end{equation}
    We now use the standard identity
    \[(1-4x)^{-\frac{1}{2}}=\sum_{n\geq0}\binom{-\frac{1}{2}}{n}(-4x)^n=\sum_{n\geq0}\binom{2n}{n}x^n.\]
    Substituting this expansion into Equation~\ref{eq:genFunctionTop} and collecting terms, we finally get
    \[T_1(x)=\sum_{n\geq3}\binom{2n-3}{n-3}x^n\]
    which proves the theorem.
\end{proof}

For a list of other combinatorial interpretations of sequence~\ref{eq:topTr}, see \href{https://oeis.org/A002054}{[A002054]}.

Notice that the same method can be used to retrieve the sequence $t_{n,\lambda}$ for any given $\lambda\in\mathcal{Y}_n$. The only difference is that for the given $\lambda$ one should carefully describe the set of polygons that contain $a_{n+1}$ in their interior.

\bibliographystyle{plain}
\bibliography{bibliography.bib}
\end{document}